\pdfoutput=1
\documentclass[
 a4paper,
 final, 
 disable,
]{article}

\usepackage{artmacs}

\usepackage{sagetex}
\usepackage{fancyvrb}

\newcommand{\F}{F} 

\newcommand{\ttf}[1]{^{\underline{#1}}}

\newcommand{\ffr}[2]{b_{#1, #2}}
\newcommand{\ffm}[2]{b_{#1, #2}}

\newcommand{\Pall}[1]{P_{#1}^{\text{all}}}
\newcommand{\Pone}[1]{P_{#1}}

\newcommand{\Ione}[1]{I_{#1}}

\newcommand{\Rone}[1]{R_{#1}}

\newcommand{\Qone}[1]{Q_{#1}}

\newcommand{\Sone}[1]{S_{#1}}

\newcommand{\Eone}[1]{E_{#1}}

\newcommand{\Aone}[1]{A_{#1}}
\newcommand{\Aonep}[1]{\Aone{#1}^{+}}

\newcommand{\qvar}{\mathbf{q}}
\newcommand{\degq}{\deg_{\qvar}}

\newcommand{\Qz}{\QQ \bbracket{z}}
\newcommand{\Qq}{\QQ (\qvar)}
\newcommand{\Qqz}{\Qq \bbracket{z}}

\newcommand{\combclass}[1]{{\mathcal #1}}
\newcommand{\ccP}{\combclass{P}}
\newcommand{\ccI}{\combclass{I}}
\newcommand{\ccR}{\combclass{R}}
\newcommand{\ccS}{\combclass{S}}
\newcommand{\ccQ}{\combclass{Q}}

\newcommand{\ccA}{\combclass{A}}
\newcommand{\ccB}{\combclass{B}}
\newcommand{\ccC}{\combclass{C}}

\newcommand{\ngenfun}[1]{{\mathrm #1}} 
\newcommand{\nP}{\ngenfun{P}}
\newcommand{\nI}{\ngenfun{I}}
\newcommand{\nR}{\ngenfun{R}}
\newcommand{\nS}{\ngenfun{S}}
\newcommand{\nQ}{\ngenfun{Q}}
\newcommand{\nAp}{\ngenfun{A}^{+}}

\newcommand{\nA}{\ngenfun{A}}
\newcommand{\nE}{\ngenfun{E}}
\newcommand{\nB}{\ngenfun{B}}
\newcommand{\nC}{\ngenfun{C}}

\newcommand{\sgenfun}[1]{{\mathsf #1}} 
\newcommand{\sP}{\sgenfun{P}}
\newcommand{\sF}{\sgenfun{F}}
\newcommand{\sI}{\sgenfun{I}}
\newcommand{\sR}{\sgenfun{R}}
\newcommand{\sS}{\sgenfun{S}}
\newcommand{\sQ}{\sgenfun{Q}}
\newcommand{\sE}{\sgenfun{E}}
\newcommand{\sA}{\sgenfun{A}}

\newcommand{\ncoeff}[1]{{\mathrm #1}} 
\newcommand{\ncP}{\ncoeff{P}}
\newcommand{\ncI}{\ncoeff{I}}
\newcommand{\ncR}{\ncoeff{R}}

\newcommand{\scoeff}[1]{{\mathsf #1}} 
\newcommand{\scP}{\scoeff{P}}
\newcommand{\scI}{\scoeff{I}}
\newcommand{\scR}{\scoeff{R}}
\newcommand{\scA}{\scoeff{A}}
\newcommand{\scE}{\scoeff{E}}
\newcommand{\scS}{\scoeff{S}}
\newcommand{\scQ}{\scoeff{Q}}

\newcommand{\MSet}{{\mathcal{MSET}}}
\newcommand{\Seq}{{\mathcal{SEQ}}}

\allowdisplaybreaks[3]

\begin{document}

\title{Counting reducible, powerful, and relatively irreducible multivariate polynomials \\ over finite fields}
\author{Joachim von~zur~Gathen\footnotemark[2]\
\and Alfredo Viola\footnotemark[3]\
\and Konstantin Ziegler\footnotemark[2]}
\pdfauthor{Joachim von zur Gathen, Alfredo Viola, and Konstantin Ziegler}

\begin{sagesilent}
load('code/multi-def.sage')
\end{sagesilent}

\maketitle

\renewcommand{\thefootnote}{} 
\footnotetext{An Extended Abstract of this paper appeared as
  \cite{gatvio10} and the full version is to appear in \emph{SIAM Journal of
  Discrete Mathematics}.}

\renewcommand{\thefootnote}{\fnsymbol{footnote}}
\footnotetext[2]{B-IT, Universit{\"{a}}t Bonn,  D-53113 Bonn, Germany, \email{{gathen,zieglerk}@bit.uni-bonn.de}}
\footnotetext[3]{Instituto de Computaci\'on, Universidad de la
  Rep\'ublica, Montevideo, Uruguay, \email{viola@fing.edu.uy}}

\renewcommand{\thefootnote}{\arabic{footnote}}

\begin{center}
\emph{Dedicated to the memory of Philippe Flajolet}
\end{center}

\begin{abstract} We present counting methods for some special classes of
  multivariate polynomials over a finite field, namely the reducible ones, the
  $s$-powerful ones (divisible by the $s$th power of a nonconstant
  polynomial), and the relatively irreducible ones (irreducible but reducible
  over an extension field).  One approach employs generating functions,
  another one uses a combinatorial method.  They yield exact formulas and approximations with relative
  errors that essentially decrease exponentially in the input size.
\end{abstract}




\section{Introduction}

Most integers are composite and most univariate
polynomials over a finite field are reducible.  The classical results
of the Prime Number Theorem and a theorem of Gau\ss\ present
approximations saying that
randomly chosen integers up to $x$ or polynomials of degree
up to $n$ are prime or irreducible with probability about $1/\ln x$
or $1/n$, respectively.

Concerning special classes of univariate polynomials over a finite
field, \cite{zsi94} counts those with a given number of distinct roots
or without irreducible factors of a given degree.  In the same
situation, \cite{art24} counts the irreducible ones in an arithmetic
progression and \cite{hay65} generalizes
these results.  \cite{coh69a} and \cite{car87}
count polynomials with certain factorization patterns and \cite{wil69}
those with irreducible factors of given degree.  Polynomials that
occur as a norm in field extensions are studied by \cite{goglut81}.

In two or more variables, the situation changes dramatically. Most
multivariate polynomials are irreducible.  \cite{car63a} provides the
first count of irreducible multivariate polynomials.  In \cite{car65},
he goes on to study the fraction of irreducibles when bounds on the
degrees in each variable are prescribed; see also \cite{coh68}.  In
this paper, we opt for bounding the total degree because it has the
charm of being invariant under invertible linear transformations.
\cite{gaolau02} consider our problem in yet another model, namely
where one variable occurs with maximal degree.  The natural generating
function (or zeta function) for the irreducible polynomials in two or
more variables does not converge anywhere outside of the
origin. \cite{wan92} notes that this explains the lack of a simple
combinatorial formula for the number of irreducible polynomials. But
he gives a $p$-adic formula, and also a (somewhat complicated)
combinatorial formula.  For further references, see
\citet[Section~3.6]{mulpan13}.

In the bivariate case, \cite{gat08-incl-gat07}
proves precise approximations with an exponentially decreasing relative
error.  \cite{bod08} gives a recursive formula for the number of irreducible
bivariate polynomials and remarks on a generalization for more than
two variables; he follows up with \cite{bod10}.  Some further types of multivariate polynomials are examined from
a counting perspective: decomposable ones (\cite{gat10a}, \cite{boddeb09}), singular
ones (\cite{gat08-incl-gat07}), and pairs of coprime polynomials (\cite{houmul09}).

This paper provides exact formulas for the numbers of reducible, $s$-powerful,
and relatively irreducible polynomials.  The latter also yields
the number of absolutely reducible polynomials.  Of these, only reducible polynomials have been treated in
the literature, usually with much larger error terms.  The
formulas yield simple, yet precise, approximations to these numbers, with rapidly decaying
relative errors.

We use two different methodologies to obtain such bounds: generating
functions and combinatorial counting.  The usual approach, see
\cite{flased09}, of analytic combinatorics on series with integer
coefficients leads, in our case, to power series that diverge
everywhere (except at $0$).  We have not found a way to make this
work.  Instead, we use power series with symbolic coefficients, namely
rational functions in a variable representing the field size.
Several useful relations from standard analytic combinatorics carry
over to this new scenario.  In a first step, this yields in a
straightforward manner exact formulas for the numbers under
consideration (Theorems~\ref{pro:R_exact_by_recursion},
\ref{pro:Q_exact_by_recursion}, and \ref{pro:E_exact_by_recursion}).  These formulas are, however, not very
transparent.  Even the leading term is not immediately visible.

In a second step, coefficient comparisons yield easy-to-use
approximations to our numbers (Theorems~\ref{thm:R-from-gen},
\ref{thm:Q_by_gen}, and \ref{thm:E_by_gen}).  The relative error is
exponentially decreasing in the bit size of the data.  As an example,
\autoref{thm:R-from-gen} gives a ``third order'' approximation for the number of reducible
polynomials, and thus a ``fourth order'' approximation for the
irreducible ones.  The error term is in the big-Oh form and thus
contains an unspecified constant.

In a third step, a different method, namely some combinatorial
counting, yields ``second order'' approximations with explicit
constants in the error term (Theorems~\ref{thm:red},
\ref{thm:Q_by_map}, and \ref{thm:E_complete}).

Geometrically, a single polynomial corresponds to a hypersurface, that
is, to a cycle in affine or projective space, of codimension 1.  This
correspondence preserves the respective notions of reducibility.  Thus,
Sections \ref{sec:gen} and \ref{sec:red} can also be viewed as counting reducible
hypersurfaces, in particular, planar curves, and
\autoref{sec:powerful} those with an $s$-fold component.  Reducible curves
embedded in higher-dimensional spaces, parametrized by the appropriate
Chow variety, are counted in \cite{cesgat13}.

\section{Notation}
\label{sec:exact}

We work in the polynomial ring $\F [x_{1}, \dots ,x_{r}]$  in
$r\geq 1$ variables over a field $\F$ and consider polynomials
with total degree equal to some nonnegative integer $n$:
\[
\Pall{r,n}(\F) = \{ f \in \F [x_{1},\dots,x_{r}] \colon \deg f = n\}.
\]
The polynomials of degree at most $n$
form an $\F$-vector space of dimension
\[
\ffr{r}{n} =  \binom{r+n}{r}=\frac{(r+n)\ttf{r}}{r!},
\]
where the \emph{falling factorial} or \emph{Pochhammer symbol} is
\begin{equation}
  \label{eq:fall_fac}
  (r+x)\ttf{r}=(r+x) \cdot (r-1+x) \cdots (1+x),
\end{equation}
for any real $x$ and any nonnegative integer $r$, see e.g. \cite{knu92a}.  Over a finite field $\FF_q$ with
$q$ elements, we have
\[
\# \Pall{r,n}(\FF_q)=
q^{\ffr{r}{n}}-q^{\ffr{r}{n-1}}=q^{\ffr{r}{n}}(1-q^{-\ffm{r-1}{n}}).
\]

The property of a certain polynomial to be reducible, squareful or
relatively irreducible is shared with all polynomials associated to
the given one. For counting them, it is sufficient to
take one representative. We choose an arbitrary monomial order, say, the degree-lexicographic one,
so that the monic polynomials are those with leading coefficient 1, and write
\[
\Pone{r,n}(\F) = \{f \in \Pall{r,n}(\F) \colon f \text{ is monic}\}.
\]
Then
\begin{equation}
\label{eq:all}
\# \Pone{r,n} (\FF_q)= \frac{\# \Pall{r,n}(\FF_q)}{q-1} =
q^{\ffr{r}{n}-1}\frac{1-q^{-\ffm{r-1}{n}}}{1-q^{-1}}.
\end{equation}
The product of two monic polynomials is again monic.

Our exact formulas are derived using a generating series, the standard tool in analytic combinatorics as presented in
\cite{flased09} by two experts who created large parts of the theory. We first recall a few general primitives from this theory that
enable one to set up symbolic equations for generating functions
starting from combinatorial specifications.  A countable set
$\combclass{C}$ with a ``size'' function $\abs{\cdot} \colon
\combclass{C} \to \ZZ_{\geq 0}$ is called a \emph{combinatorial class}
if the preimage of any $n \in \ZZ_{\geq 0}$ is finite.  The number of elements of size $n$ is denoted by $\ncoeff{C}_{n}$ and these numbers are encoded in the \emph{generating function}
$\ngenfun{C}(z)$ of the sequence $\ncoeff{C}_{n}$:
\begin{equation}
  \ngenfun{C}(z) = \sum_{n \geq 0} \ncoeff{C}_{n} z^{n} \in \ZZ_{\geq 0} \bbracket{z}.
\end{equation}
We sometimes omit the argument $z$.  Before we tackle the task of counting polynomials, let us recall some basics about power series.  An element in the ring of univariate power series over a ring is invertible if
and only if its constant term is invertible. We call a power series
\emph{original} if its constant term vanishes, so that its graph
passes through the origin.  The power series
\begin{equation}
  \log (1-z) = - \sum_{n \geq 1} \frac{z^{n}}{n} \in \Qz \label{eq:91}
\end{equation}
is original and substituting a power series $f$ in another
power series $g$ is well-defined if $f$ is original.

Two combinatorial classes $\ccA$ and $\ccB$ are \emph{isomorphic} if there is a size-preserving bijection $\ccA \to \ccB$ or equivalently if $\nA = \nB$.  We recall three basic constructions of new combinatorial classes from given ones; see \citet[Section I. 2.]{flased09}.

Let $\ccA$
and $\ccB$ be two combinatorial classes.  We define the \emph{disjoint union}
\begin{equation}
  \label{eq:38}
  \ccA \, \dot{\cup}\, \ccB = \{ \{0\} \times \ccA \} \cup \{ \{1\} \times \ccB\}.
\end{equation}
The size of an element $(0,a)$ or $(1,b)$ is defined as the size of $a$ or $b$, respectively.  We also define the \emph{sequence class}
\begin{equation}
  \label{eq:39}
  \Seq (\ccA) = \{ (\alpha_{1}, \dots, \alpha_{\ell}) \colon \ell \geq 0, \alpha_{i} \in \ccA\},
\end{equation}
where $\abs{(\alpha_{1}, \dots, \alpha_{\ell})} = \sum_{i} \abs{\alpha_{i}}$.  This is a combinatorial class, if $\ccA$ contains no element of size $0$.  Finally, we derive the \emph{multiset class}
\begin{equation}
  \label{eq:40}
  \MSet (\ccA) = \Seq (\ccA) /\mathbin{\sim},
\end{equation}
where $(\alpha_{1}, \dots, \alpha_{\ell}) \sim (\beta_{1}, \dots, \beta_{\ell})$ if there is a permutation $\sigma$ of $\{1, \dots, \ell\}$ such that $\alpha_{i} = \beta_{\sigma (i)}$ for all $i$.  This class contains all finite sequences of elements from $\ccA$ where repetition is allowed, but ordering ignored.  The generating functions for these constructions are
classic applications of combinatorics.
\begin{fact}[see {\citet[Theorems~I.1 and I.5]{flased09}}]
\label{fac:gf}
Let $\ccA$, $\ccB$, and $\ccC$  be combinatorial classes.
  \begin{ronumerate}
  \item\label{item:1} If $\ccA = \ccB \, \dot{\cup}\,  \ccC $, then $\nA = \nB + \nC$.
  \item\label{item:2} If $\ccA = \MSet (\ccB)$ and $\ncoeff{B}_{0} =0$, then
    \begin{equation}
      \label{eq:41}
      \ngenfun{B} = \sum_{k \geq 1} \frac{\mu (k)}{k} \log (\ngenfun{A}(z^{k})),
    \end{equation}
where $\mu$ is the number-theoretic M\"{o}bius-function, defined as
\begin{equation}
  \label{eq:81}
\mu(k) = \begin{cases*}
1 & if $k = 1$, \\
(-1)^{\ell} & if $k$ is the product of $\ell$ distinct primes, \\
0 & otherwise.
\end{cases*}
\end{equation}
\end{ronumerate}
\end{fact}

\section{Generating functions for reducible polynomials}
\label{sec:gen}

To study reducible polynomials, we consider the following subsets of $\Pone{r,n}(\F)$:
\begin{align*}
  \Ione{r,n}(\F) & = \{ f \in \Pone{r,n}(\F) \colon f \textrm{ irreducible} \}, \\
  \Rone{r,n}(\F) & = \Pone{r,n}(\F) \setminus \Ione{r,n}(\F).
\end{align*}
In the usual notions, the polynomial $1$ is neither reducible nor
irreducible.  In our context, it is natural to have $\Rone{r,0}(\F) =
\{ 1 \}$ and $ \Ione{r,0}(\F) = \varnothing$.

The sets of polynomials
\begin{align}
  \ccP & = \bigcup_{n \geq 0} \Pone{r,n} (\FF_{q}), \\
  \ccI & = \bigcup_{n \geq 0} \Ione{r,n} (\FF_{q}), \\
  \ccR & = \ccP \setminus \ccI,
\end{align}
are combinatorial classes with the total degree as size functions and
we denote the corresponding generating functions by $\nP, \nI, \nR \in
\ZZ_{\geq 0} \bbracket{z}$, respectively.  Their coefficients are
\begin{align}
\nP_{n} & = \nP_{r,n} (\FF_{q}) = \# \Pone{r,n} (\FF_{q}) = q^{\ffr{r}{n}-1} \frac{1-q^{-\ffm{r-1}{n}}}{1-q^{-1}}, \label{eq:44} \\
\nR_{n} & = \nR_{r,n} (\FF_{q}) = \#  \Rone{r,n} (\FF_{q}), \\
\nI_{n} & = \nI_{r,n} (\FF_{q}) = \#  \Ione{r,n} (\FF_{q}), \\
\end{align}
respectively, dropping $\FF_{q}$ and $r$ from the notation.
By definition, $\ccP$ is isomorphic to the disjoint union of $\ccR$ and $\ccI$, and therefore
\begin{equation}
\label{eq:42}
\nR = \nP - \nI
\end{equation}
by \autoref{fac:gf} \ref{item:1}.  By unique factorization, every element in $\ccP$ corresponds to an unordered finite sequence of irreducible polynomials, where repetition is allowed.  Hence $\ccP$ is isomorphic to $\MSet (\ccI)$ and  by \autoref{fac:gf} \ref{item:2},
\begin{equation}
\nI  = \sum_{k \geq 1} \frac{\mu (k)}{k} \log \nP (z^{k}). \label{eq:28}
\end{equation}

A Maple implementation of the resulting algorithm to compute the number of reducible polynomials is described in
\autoref{fig:algoirr}.  It is easy to program and execute and was used to
calculate the number of bivariate reducible polynomials in
\citet[Table~2.1]{gat08-incl-gat07}.  We extend these exact results in \autoref{tab:R}.

\begin{figure}
\centering
\begin{Verbatim}[frame=single]
allpolysGF:=proc(z,N,r) local i: option remember:
    sum('simplify((q^binomial(r+i,r)-q^binomial(r+i-1,r))/
        (q-1))*z^i',i = 0..N):
end:

irreduciblesGF:=proc(z,N,r) local k: option remember:
    convert(taylor((sum('mobius(k)/k*log(allpolysGF(z^k,N,r))',
        k=1..N)), z, N+1), polynom):
end:

reduciblesGF:=proc(z,N,r) option remember:
    allpolysGF(z,N,r)-irreduciblesGF(z,N,r):
end:

reducibles:=proc(n,r)
coeff(sort(expand(reduciblesGF(z,n,r))),z^n):
end:
\end{Verbatim}
\caption{Maple program to compute the number of
monic reducible polynomials in $r$ variables of degree $n$.}
\label{fig:algoirr}
\end{figure}

\begin{table}\footnotesize
\centering
\sagestr{R_table(3,6)}

\sagestr{R_table(4,4)}

\sagestr{R_table(5,4)}
  \caption{Exact values of $\# \Rone{r,n} (\FF_{q})$ for small values
of $r$ and $n$.  For $n<4$, these are the numbers given in \autoref{thm:R-from-gen}.}
\label{tab:R}
\end{table}

This approach quickly leads to explicit formulas.  For a positive
integer $n$, a \emph{composition} of $n$  is a sequence $j = (j_{1},
j_{2}, \dots, j_{\abs{j}})$ of positive integers $j_{1}, j_{2}, \dots,
j_{\abs{j}}$ with $j_{1} +
j_{2} + \dots + j_{\abs{j}} = n$, where $\abs{j}$ denotes the length
of the sequence.  We define the set
\begin{equation}
  \label{eq:10}
M_{n}  = \{ \text{compositions of $n$} \}.
\end{equation}
  This standard
combinatorial notion is not to be confused with the composition of
polynomials, for which also counting results are available.

\begin{theorem}
\label{pro:R_exact_by_recursion} For $r \geq 1$, $q \geq 2$, and $\nP_{n}$ as in \eqref{eq:44}, we have
\begin{equation}
  \label{eq:12}
  \begin{split}
    \nI_{0} & = 0, \\
    \nI_{n} & = - \sum_{k \,\mid\, n} \frac{\mu (k)}{k} \sum_{j \in
      M_{n/k}} \frac{(-1)^{\abs{j}}}{\abs{j}} \nP_{j_{1}} \nP_{j_{2}} \cdots \nP_{j_{\abs{j}}},
  \end{split}
\end{equation}
for $n \geq 1$, and therefore
  \begin{align}
    \nR_{0} & = 1, \\
    \nR_{n} & = \nP_{n} + \sum_{k \,\mid\, n} \frac{\mu (k)}{k} \sum_{j \in M_{n/k}} \frac{(-1)^{\abs{j}}}{\abs{j}} \nP_{j_{1}} \nP_{j_{2}} \cdots \nP_{j_{\abs{j}}},
  \end{align}
for $n \geq 1$.
\end{theorem}

\begin{proof}
We consider the original power series $\ngenfun{F} = 1 - \sP = - \sum_{i \geq 1} \nP_{i} z^{i}$.  The Taylor expansion \eqref{eq:91} of
$\log(1-\ngenfun{F}(z^k))$ in \eqref{eq:28} yields
\begin{align}
\ngenfun{I} &= - \sum_{k \geq 1} \frac{\mu(k)}{k} \sum_{i \geq 1} \frac{\ngenfun{F}(z^k)^i}{i} = - \sum_{k \geq 1} \frac{\mu(k)}{k} \sum_{i \geq 1} \frac{(-1)^{i}}{i} \Bigl(\sum_{j \geq 1} \nP_{j} z^{jk}\Bigr)^i \\
&= - \sum_{k \geq 1} \frac{\mu(k)}{k} \sum_{i \geq 1} \frac{(-1)^{i}}{i} (\nP_{1}  z^{k} + \nP_{2}  z^{2k} + \nP_{3}  z^{3k} + \dots)^i, \\
\nI_{0} & = 0, \\
\nI_{n} & = - \sum_{k \,\mid\, n} \frac{\mu(k)}{k} \sum_{i \geq 1} \frac{(-1)^{i}}{i} \sum_{\substack{j \in M_{n/k} \\ \abs{j}=i}} \nP_{j_{1}} \nP_{j_{2}} \cdots \nP_{j_{i}},
\end{align}
for $n \geq 1$, which proves the claimed formulas for $\nI$.  The
results for $\nR$ follow by \eqref{eq:42}.
\end{proof}

We check that the formula yields the well-known one, see \citet[Theorem 3.25]{lidnie97}, in the univariate case, where $r=1$.  We then have $\nP_j = q^j$ and so
$\nP_{j_{1}} \nP_{j_{2}} \cdots \nP_{j_{i}} = q^{n/k}$ for any composition $j_{1} + j_{2} +  \dots + j_{i}
= n/k$.  Moreover, the number of compositions of $m$ with $i$ components is $\binom{m-1}{i-1}$, see \citet[Section I.3.1]{flased09}.  As a consequence we have for $k$ dividing $n$
\begin{align}
\sum_{\substack{j \in M_{n/k} \\ \abs{j}=i}}
\frac{(-1)^{i}}{i} \nP_{j_{1}} \nP_{j_{2}} \cdots \nP_{j_{i}} & = q^{n/k} \sum_{i \geq 1} \frac{(-1)^{i}}{i} \binom{n/k-1}{i-1} \\
& = \frac{ k q^{n/k}}{n} \sum_{i \geq
1} (-1)^i \binom{n/k}{i} = - \frac{k q^{n/k}}{n}, \\
\nI_{n}  & = \frac{1}{n}\sum_{k \,\mid\, n} \mu(k) q^{n/k}. \label{eq:80}
\end{align}

\cite{coh68} notes that, compared to the univariate case,
``the situation is different and much more difficult. In this case, no
explicit formula [...] is available.''

For $r \geq 2$, the power series $\nP$, $\nI$, and $\nR$ do not converge anywhere except at 0, and the standard asymptotic arguments of analytic combinatorics are inapplicable.  We now deviate from this approach and move from power series in $\Qz$ to power series in $\Qqz$, where $\qvar$ is a symbolic variable representing the field size.  For $r \geq 2$ and $n \geq 0$ we let
\begin{equation}
  \label{eq:19}
\scP_{n}(\qvar) = \scP_{r,n}(\qvar) = \qvar^{\ffr{r}{n}-1}\frac{1-\qvar^{-\ffm{r-1}{n}}}{1-\qvar^{-1}} \in \ZZ [\qvar],
\end{equation}
where we usually omit $r$ from the notation.  As examples, we have
\begin{equation}
  \label{eq:29}
\scP_{0} (\qvar) = 1,   \scP_{1} (\qvar) = \qvar^{r} \frac{1 - \qvar^{-r}}{1-\qvar^{-1}}, \text{ and } \scP_{2}(\qvar)  = \qvar^{r(r+3)/2} \frac{1- \qvar^{-r(r+1)/2}}{1-\qvar^{-1}}.
\end{equation}
We define the power series $\sP, \sI, \sR \in \Qqz$ by
\begin{align}
  \sP (\qvar, z) & = \sum_{n \geq 0} \scP_{n}(\qvar) z^{n}, \label{eq:5} \\
  \sI (\qvar, z) & = \sum_{k \geq 1} \frac{\mu (k)}{k} \log \sP (\qvar, z^{k}), \label{eq:relationIP} \\
  \sR (\qvar, z) & = \sP(\qvar, z) - \sI(\qvar, z). \label{eq:70}
\end{align}
Now $1 - \sgenfun{P}(\qvar, z^{k})$ is an
original power series, and $\log \sgenfun{P}(\qvar, z^{k})$ and  $\sI$
are well-defined, with $\sI (\qvar, 0)=0$.  For $q \in \QQ$, the
rational functions in $\QQ (\qvar)$ without pole at $\qvar \gets q$
form a ring, the localization $\QQ[\qvar]_{(\qvar-q)}$.  If we restrict the power series coefficients to this ring, the evaluation map which substitutes an integer $q$ for $\qvar$ is a ring
homomorphism. Since $\scP_{n}$ is
actually a polynomial in $\qvar$, this poses no restriction in our case, and
evaluating $\qvar \gets q$ maps $\sP(\qvar, z)$ to $\nP(z)$
coefficientwise.  In other words, the coefficient of $z^{n}$ equals
\begin{equation}
  \label{eq:20}
  [z^{n}] \sP (q,z) = \ncP_{n}
\end{equation}
by \eqref{eq:all}.  Furthermore, $\sI$ and $\sR$ relate to $\sP$ in the same way as $\nI$ and $\nR$ do to $\nP$, so that
\begin{align}
  [z^{n}] \sI (q,z) = \ncI_{n}, \\
  [z^{n}] \sR (q,z) = \ncR_{n}.
\end{align}

The formula of \autoref{pro:R_exact_by_recursion} is exact but somewhat cumbersome.  A main goal in this work is to find simple yet precise approximations, with rapidly decaying error terms.   We fix some notation.
For nonzero $f \in \QQ(\qvar)$, $\degq f$ is the degree of $f$, that is, the numerator
degree minus the denominator degree.  Thus $\degq \scP_{n} = b_{r,n}-1$ and $\degq (f+g) \leq \max \{ \degq f, \degq g\}$.
The appearance of $O(\qvar^{-m})$ with a positive integer $m$ in an
equation means the existence of some $f$ with degree at most $-m$ that makes the
equation valid.  The charm of our approach is that we obtain results for any
``fixed'' $r$ and $n$.  If a term $O(\qvar^{-m})$ appears, then we may
conclude a numerical
asymptotic result for growing prime powers $q$.

We start with a degree comparison for certain products of the
$\scP_{i} (\qvar)$ and sometimes omit the argument $\qvar$.

\begin{lemma}
Let $r \geq 2$ and $n \geq 0$.
\label{lem:deg_q}
  \begin{ronumerate}
  \item \label{item:degi} For $i,j \geq 0$, we have $\degq (\scP_i \cdot
    \scP_j) \leq \degq \scP_{i+j}$, with equality if and only if  $ij = 0$.
  \item  \label{item:degii} For $1 \leq k \leq n/2$, the sequence of integers
    $\degq (\scP_k \cdot \scP_{n-k})$ is strictly decreasing in $k$.
  \item \label{item:degiii} For $3 \leq k \leq n/2$, we have $\degq \scP_{1}^{2}\scP_{n-2} \geq \degq \scP_{k} \scP_{n-k}$, with equality only for $(r,n, k) = (2,6, 3)$.
  \end{ronumerate}
\end{lemma}

\begin{proof}
  \begin{ronumerate}
  \item The claimed inequality is equivalent to
    \begin{equation}
      \label{eq:1}
      \binom{r+i}{r} + \binom{r+j}{r} - 1 \leq \binom{r+i+j}{r},
    \end{equation}
    which follows by considering the choices of $r$-element subsets
    from a set with $r+i+j$ elements.  Since $r\geq 2$, this inequality is
    strict if and only if both $i$ and $j$ are nonzero.
  \item Using \eqref{eq:19}, we define a function $u$ as
    \begin{equation}
      \label{eq:2}
u(k) = \degq (\scP_k \cdot \scP_{n-k}) = \binom{r+k}{r} + \binom{r+n-k}{r}-2.
    \end{equation}
    We extend the domain of $u(k)$ to real
    numbers $k$ between $1$ and $n/2$ by means of falling factorials
    as in \eqref{eq:fall_fac}
\[
u(k) = \frac{(r+k)\ttf{r}}{r!} + \frac{(r+n-k)\ttf{r}}{r!} - 2.
\]
  It is sufficient to show that the affine transformation $\bar{u}$
  with
  \begin{equation*}
    \bar{u}(k) =r!\cdot (u (k)+2)  =(r+k)\ttf{r}+(r+n-k)\ttf{r}
  \end{equation*}
  is strictly decreasing.  The first derivative with respect to $k$ is
  \begin{equation}
    {\bar{u}}'(k) = \sum_{1\leq i \leq r} \Bigl( \frac{(r+k)\ttf{r}}{i+k} - \frac{(r+n-k)\ttf{r}}{i+n-k} \Bigr).
  \end{equation}
  Since $0< i+k < i+n-k$ for $1 < k < n/2$, each difference is
  negative, and so is ${\bar{u}}'(k)$.
  \item  Since $r \geq 2$ we have
\begin{align}
 (r-2)(r-1)(r+5) & \geq 0, \\
b_{r-1,4} & \geq b_{r,3} -2r -1, \label{eq:18} \\
 2r + b_{r,4} -1 & \geq 2b_{r,3}-2, \\
\degq \scP_{1}^{2}\scP_{4} & \geq \degq \scP_{3}\scP_{3},
\end{align}
and equality if and only if $r=2$.  This proves the claimed inequality for $n=6$.

For $n>6$ we have $b_{r-1,n-2} > b_{r-1,4}$ and with \eqref{eq:18} follows
\begin{align}
  b_{r-1,n-2} & > b_{r,3} - 2r -1, \\
  2r+b_{r,n-2} -1 & > b_{r,3} + b_{r,n-3} -2, \\
\degq \scP_{1}^{2}\scP_{n-2} & > \degq \scP_{3} \scP_{n-3},
\end{align}
which proves \ref{item:degiii} for $k=3$ and by the monotonicity
proven in \ref{item:degii} also for all larger $k$. \qedhere
  \end{ronumerate}
\end{proof}

\begin{theorem}\label{thm:R-from-gen}
  Let $r\geq 2$ and
  \begin{equation}
    \rho_{r,n}(\qvar)
= \qvar^{\binom{r+n-1}{r}+r-1}
    \frac{1-\qvar^{-r}}{(1-\qvar^{-1})^2}  \in \QQ(\qvar). \label{eq:rho}
\end{equation}
Then
\begin{align}
\scoeff{R_{0}} & = 1, \\
\scoeff{R_{1}} & = 0, \\
\scoeff{R_2} & = \frac{ \rho_{r,2}(\qvar) }{2} \cdot (1-\qvar^{-r-1}), \\
\scoeff{R_3} & = \rho_{r,3}(\qvar) \Bigl( 1-\qvar^{-r(r+1)/2} +
\qvar^{-r(r-1)/2}
\frac{1-2\qvar^{-r}+2\qvar^{-2r-1}-\qvar^{-2r-2}}{3(1-\qvar^{-1})}\Bigr),
\\
\scoeff{R_4} & = \rho_{r,4} (\qvar) \cdot \Bigl( 1 +
\qvar^{-\binom{r+1}{3}} \cdot  \frac{1 +
  O(\qvar^{-r(r-1)/2})}{2 (1-\qvar^{-r})} \Bigr), \label{eq:48} \\
\intertext{and for $n \geq 5$}
\scoeff{R_n} & = \rho_{r,n} (\qvar) \cdot \Bigl( 1 + \qvar^{-\binom{r+n-2}{r-1}+r(r+1)/2} \cdot  \frac{1+ O(\qvar^{-r(r-1)/2})}{1-\qvar^{-r}} \Bigr). \label{eq:26}
\end{align}
\end{theorem}

\begin{proof}
We start the symbolic analog of our approach in the proof of \autoref{pro:R_exact_by_recursion} with the original power series $\sgenfun{F}= 1 - \sP  = -\sum_{i \geq 1} \scP_{i}z^{i}$.  The Taylor expansion of $\log (1-\sgenfun{F} (z^{k}))$ in \eqref{eq:relationIP} yields
\begin{equation}
 \label{eq:21}
\sR = \sP - \sI = 1 +  \sum_{i \geq 2} \frac{\sgenfun{F}^i}{i} + \sum_{k \geq 2} \frac{\mu(k)}{k}
\sum_{i \geq 1} \frac{\sgenfun{F}(z^k)^i}{i}.
\end{equation}

Since $\scR_{n} = [z^n] \sR$, we find $\scR_{0}=1$, $\scR_{1}=0$,
$\scR_{2} = (\scP_{1}^2+\scP_{1})/2$, and $\scR_{3}= \scP_{2}\scP_{1}
- (\scP_{1}^3-\scP_{1})/3$.  Together with \eqref{eq:29}, these imply
the claims for $n < 4$.

\begin{table}[h!]\footnotesize
\centering
\begin{tabular}{|L|L|L|L|}
  \hline
  & i & \text{summands} & \text{summands with} \\
& & & \text{largest degree in } \qvar \\
\hline
  [z^n] \sF^i & 2 & \scP_j \scP_{n-j}, & \scP_1\scP_{n-1}, \scP_2\scP_{n-2},  \\
  & & \quad 1 \leq j \leq n/2 & \quad  \scP_3\scP_{n-3} \text{ (for $n \geq 6$)}\\
  & \geq 3 & \scP_{j_1} \scP_{j_2} \cdots \scP_{j_i},  & \scP_1^2\scP_{n-2} \\
  & & \quad 1 \leq j_1 \leq j_2 \leq \dots \leq j_i \leq n,  & \\
  & & \quad j_1 + j_2 + \dots + j_i = n, & \\
  {[z^n]} \sF(z^k)^i & 1 & \scP_{n/k} & \scP_{n/k} \\
  & \geq 2 & \scP_{j_1} \scP_{j_2} \cdots \scP_{j_i}, &  \scP_1\scP_{n/k-1} \\
  & & \quad 1 \leq j_1 \leq j_2 \leq \dots \leq j_i \leq n/k,  & \\
  & & \quad j_1 + j_2 \dots + j_i = n/k, & \\
  \hline
\end{tabular}
\caption{Summands of $\scR$ and bounds on their degrees in $\qvar$.}
\label{tab:terms_of_R}
\end{table}

When $n\geq 4$, the contributions to $[z^{n}] \sR$ from both sums in \eqref{eq:21} are displayed
in \autoref{tab:terms_of_R}, distinguishing the smallest possible
value for $i$ from the remaining larger ones.  The third column lists all summands.  We first show that the last column displays the terms of largest degree in their row, and then compare the summands in the last column.  The terms of $[z^{n}]\sF^{i}$ are products of $i$ factors
\[
\scP_{j_1} \scP_{j_2} \cdots \scP_{j_i}, \quad 1 \leq j_1 \leq j_2 \leq \dots \leq j_i \leq n,
\]
with $j_1 + j_2 + \dots + j_i = n$.   For $i = 2$, we find
\begin{equation}
\label{eq:23}
\degq  \scP_1\scP_{n-1} > \degq \scP_2\scP_{n-2} > \degq  \scP_j\scP_{n-j}
\end{equation}
for all $j$ with $3 \leq j \leq n/2$ by \autoref{lem:deg_q}~\ref{item:degii}. For $i \geq 3$,
\begin{equation}
\degq \scP_{1}^{2}\scP_{n-2} \geq \degq \scP_{j_1} \scP_{j_2} \cdots \scP_{j_i}
\end{equation}
for all admissible values of $j_1, \dots, j_i$ by repeated application
of \autoref{lem:deg_q}~\ref{item:degi} and a single instance of
\ref{item:degii}.  Let $k$ divide $n$.  Then $[z^{n}]\sF (z^{k}) = -\scP_{n/k}$ and $[z^{n}]\sum_{i \geq 2} \sF (z^{k})^{i}$ has degree $\degq \scP_{1}\scP_{n/k-1}$ as shown above for $k=1$.

We continue the comparison started in \eqref{eq:23} by noting that
$\degq \scP_2\scP_{n-2} > \degq \scP_{1}^{2}\scP_{n-2}$ by
\autoref{lem:deg_q}~\ref{item:degi}, and also $\degq
\scP_{1}^{2}\scP_{n-2} \geq \degq \scP_{j} \scP_{n-j}$ for all $3 \leq j \leq n/2$ with equality only for $(r,n,j) = (2,6,3)$ by \autoref{lem:deg_q}~\ref{item:degiii}. Furthermore, since $\degq \scP_1 \geq 1$, we have for $k \geq 2$
\begin{equation}
  \label{eq:30}
\degq \scP_1^2 \scP_{n-2} > \degq \scP_{n-2} \geq   \degq \scP_{n/k} > \degq \scP_{1}\scP_{n/k-1},
\end{equation}
by \autoref{lem:deg_q} \ref{item:degi}.  Therefore, the summands of largest degree in $\qvar$ are in decreasing order $\scP_1 \scP_{n-1}$, $\scP_2 \scP_{n-2}$, and $\scP_1^2\scP_{n-2}$. For $n=4$, this leads to
\begin{align}
  \scoeff{R}_4 & = \scP_{1}\scP_{3}+\scP_{2}^2/2 -
  \scP_1^2\scP_2(1+O(\qvar^{-1})) \\
  & = \scP_1 \scP_3 \biggl(1 + \frac{\scP_2^2}{2\scP_1\scP_3}\cdot
    \Bigl(1 - \frac{\scP_1^2}{\scP_2} \cdot (1+O(\qvar^{-1})) \Bigr)
  \biggr),
\end{align}
 while for $n \geq 5$, $(r,n) \neq (2,6)$ we have
 \begin{align}
\scoeff{R_n} & = \scP_{1}\scP_{n-1}+\scP_{2}\scP_{n-2} - \scP_1^2
\scP_{n-2}(1+O(\qvar^{-1}))  \\
& =
\scP_{1}\scP_{n-1}\biggl(1+\frac{\scP_{2}\scP_{n-2}}{\scP_{1}\scP_{n-1}}\cdot
  \Bigl( 1 - \frac{\scP_1^2}{\scP_2}(1+O(\qvar^{-1})) \Bigr)
\biggr). \label{eq:108}
 \end{align}
For $(r,n)=(2,6)$, we have \eqref{eq:108} with $(1/2 + O(\qvar^{-1}))$ instead of
$(1+O(\qvar^{-1}))$.

The estimates \eqref{eq:48} and \eqref{eq:26} follow from
\begin{align}
\scP_{1}\scP_{n-1} & = \rho_{r,n} (\qvar) (1 - q^{-b_{r-1,n-1}}), \\
\frac{\scP_{2}\scP_{n-2}}{\scP_{1}\scP_{n-1}} & = q^{-b_{r-1,n-1}+b_{r-1,2}} \frac{1+O(\qvar^{-r(r-1)/2})}{1-q^{-r}}, \text{ and} \\
\frac{\scP_{1}^{2}}{\scP_{2}} & = O(\qvar^{-r(r-1)/2}). \qedhere
\end{align}
\end{proof}

\cite{ale06} lists $(\#
I_{r,n}(\Fq))_{n \geq 0}$ as A115457--A115472 in The On-Line Encyclopedia
of Integer Sequences, for $2 \leq r \leq 6$ and prime $q \leq 7$.

\citet[Theorem~7]{bod08} states (in our notation)
\[
1-\frac{\# I_{r, n}}{ \# P_{r, n} } \sim q^{-b_{r-1,n}-r}
\frac{1-q^{-r}}{1-q^{-1}}.
\]
\cite{houmul09} provide results for
$\# I_{r,n} (\FF_q)$.  These do not yield error bounds for the approximation
of $\# R_{r,n} (\FF_q)$.  \cite{bod10} also uses \eqref{eq:28}.  Without proving the required bounds on the various terms, as in \autoref{lem:deg_q}, he claims a result similar to \eqref{eq:26}, but only for values of
$n$ that tend to infinity and with an unspecified multiplicative factor $O(1)$ in the place of our $(1 + O(\qvar^{-r(r-1)/2}))/(1-\qvar^{-r})$ in the error term; the latter is independent of $n$.

Our approach can be described as follows.  We start in the usual
framework of algebraic combinatorics with a power series, $\ngenfun{P}
= \sum_{n \geq 0} \ncoeff{P_{n}} z^{n}$ in our case, with
well-known integer coefficients. Then we consider a well-defined
series, $\ngenfun{I} = \sum_{n \geq 0} \ncoeff{I_{n}}z^{n}$ in our case,
whose coefficients we want to determine. We find a description of
$\ngenfun{P}$ as $f(\ngenfun{I})$ and turn this around to get $\ngenfun{I} =
g(\ngenfun{P})$, usually by M\"obius inversion.  For convergent series,
we can then apply powerful tools from calculus, such as singularity
analysis, to analyze the asymptotic behavior of the coefficients.

Since our series are not convergent, we  deviate from the standard
approach.  The coefficients $\ncoeff{P_{n}}$ are rational functions of
the field size $q$.  We introduce a variable $\qvar$ and define a
power series $\sgenfun{P} \in \QQ (\qvar) \bbracket{z}$, whose
coefficients are rational functions in a variable $\qvar$, such that $\sgenfun{P}(q,z) =
\ngenfun{P}$.  Then $g (\sgenfun{P})$ is well-defined, and we set
$\sgenfun{I} = g(\sgenfun{P}) \in \QQ (\qvar) \bbracket{z}$.  Then
$[z^{n}] \, \sgenfun{I}(q, z) = \ncoeff{I_{n}}$.  We now estimate the
degrees of the terms in $g(\sgenfun{P})$.  This yields $\sgenfun{I} =
h(\qvar) (1+ O(\qvar^{-m}))$, with a main contribution $h(\qvar) \in \QQ
(\qvar)$ and a relative error $O(\qvar^{-m})$, which is an unspecified
rational function of degree at most $-m$.

Overall, we first have to determine $\ngenfun{P}, \ngenfun{I}, f$, and
$g$, which is often a substantial  part of the labor in the standard
framework.  From then on, our derivation enjoys three advantages.
\begin{itemize}
\item No convergence of the power series is required.
\item A clean concentration on the degrees of the various
  contributions, as embodied in Lemmas~\ref{lem:deg_q}, \ref{lem:conv}, and \ref{lem:w_mon}.
\item The degree of a sum of rational functions is bounded by the
  degree of the summands.
\end{itemize}
In the standard approach, the bound for a sum as in the third point
has to be multiplied by the number of summands.  As to the second
point, one sometimes sees in the literature
a simple claim of what the main contribution is, without
argument.  It is not clear whether this constitutes a mathematical
proof in the usual sense.  Since our series are not convergent, the first point is a definitive requirement.

\section{Explicit bounds for reducible polynomials}
\label{sec:red}

We now describe a third approach to counting the reducible
polynomials.  The derivation is somewhat more involved.  The payoff of this
additional effort is an explicit relative error bound in \autoref{thm:red}.
However, the calculations are sufficiently complicated for us to stop at the
first error term.  Thus we replace the asymptotic $1 +
  O(\qvar^{-r(r-1)/2})$ in \autoref{thm:R-from-gen} by $1/(1-q^{-1})$.

We consider, for integers $1 \leq k < n$, the sets
\begin{equation}
  \Rone{r,n,k} (F)= \{g \cdot h \colon g \in \Pone{r,k}(F), h \in \Pone{r,n-k} (F) \} \subseteq \Pone{r,n} (F).
\end{equation}
For the remainder of this section we restrict ourselves to finite fields $\FF_{q}$, which we omit from the notation.  Then
\begin{equation}
\label{ineq:size_of_im}
  \# \Rone{r,n,k}  \leq \# \Pone{r,k} \cdot
  \#\Pone{r,n-k}  = q^{u(k)}
  \frac{(1-q^{-\ffr{r-1}{k}})(1-q^{-\ffr{r-1}{n-k}})}{(1-q^{-1})^2},
\end{equation}
with $u(k)=b_{r,k}+b_{r,n-k}-2$ as in \eqref{eq:2}.  The
asymptotic behavior of this upper bound is dominated by the behavior
of $u(k)$.  Since $\Rone{r,n,k} = \Rone{r,n,n-k}$, we assume without
loss of generality $k \leq n/2$.  From \autoref{lem:deg_q} \ref{item:degii}, we know that, for any
$r, n\geq 2$, $u(k)$ is strictly decreasing for $1
\leq k \leq n/2$.  As $u(k)$ takes only integral values for
integers $k$ we conclude that
\begin{equation}
  \label{eq:geosum}
  \sum_{2\leq k \leq n/2} q^{u(k)} < q^{u(2)} \sum_{k\geq 0}
  q^{-k} = \frac{q^{u(2)}}{1-q^{-1}}.
\end{equation}

\begin{theorem}
  \label{thm:red}
Let $r, q \geq 2$, and $\rho_{r,n}$ as in
  \autoref{thm:R-from-gen}. We have
  \begin{align}
  \# \Rone{r,0} (\FF_{q}) & = 1, \\
   \# \Rone{r,1} (\FF_{q}) & = 0, \\
    \# \Rone{r,2} (\FF_q) & = \frac{ \rho_{r,2}(q)}{2} \cdot (1-q^{-r-1}), \\
    \left| \# \Rone{r,3} (\FF_q)- \rho_{r,3}(q) \right| &  = \rho_{r,3}(q) \cdot q^{-r(r-1)/2}
      \frac{1-2q^{-r}+2q^{-2r-1}-q^{-2r-2}}{3(1-q^{-1})} \label{eq:27} \\
& \leq
    \rho_{r,3}(q) \cdot q^{-r(r-1)/2}, \\
\intertext{and for $n \geq 4$}
  \left| \# \Rone{r,n} (\FF_q)-
      \rho_{r,n}(q) \right| & \leq \rho_{r,n}(q) \cdot
    \frac{q^{-\binom{r+n-2}{r-1}+r(r+1)/2}}{(1-q^{-1})(1-q^{-r})} \label{eq:16} \\
& \leq \rho_{r,n}(q) \cdot 3 q^{-\binom{r+n-2}{r-1}+r(r+1)/2}.
\end{align}
\end{theorem}

\begin{proof}
For $n<4$, the claims follow from \autoref{thm:R-from-gen}. We remark that the fraction on the right-hand side of \eqref{eq:27} is actually bounded by $2/3$.  For $n \geq 4$, the proof proceeds in three steps.  We claim
\begin{align}
      \# \Rone{r,n} & \leq \rho_{r,n}(q)
      \Bigl(1+ \frac{q^{-\ffr{r-1}{n-1}+\ffr{r-1}{2}}}{(1-q^{-1})(1-q^{-r})}\Bigr), \label{eq:33}\\
      \#\Ione{r,n} & \geq \# \Pone{r,n} \Bigl( 1-3 q^{-b_{r-1,n}+r} \frac{1-q^{-r}}{1-q^{-1}}\Bigr), \label{eq:34}\\
      \# \Rone{r,n} & \geq \rho_{r,n}(q)
      \Bigl(1-3q^{-b_{r-1,n-1}+r}\frac{1-q^{-r-1}}{1-q^{-1}}\Bigr) \label{eq:35}.
\end{align}

We start with the proof of \eqref{eq:33}. Using $\Rone{r,n} = \bigcup_{1\leq k\leq n/2} \Rone{r,n,k}$ and inequality~\eqref{ineq:size_of_im}, we have
  \begin{align}
    \# \Rone{r,n} & \leq \sum_{1 \leq k \leq n/2} \# \Rone{r,n,k}  \\
    & \leq \frac{1}{(1-q^{-1})^2} \sum_{1\leq k \leq n/2}
    q^{u(k)} (1-q^{-\ffr{r-1}{k}})(1-q^{-\ffr{r-1}{n-k}}) \\
    & < \frac{1}{(1-q^{-1})^2} \sum_{1\leq k \leq n/2}
    q^{u(k)} (1-q^{-\ffr{r-1}{k}}). \label{eq:R_up}
  \end{align}
  For the sum, \eqref{eq:geosum} shows
  \begin{align}
\sum_{1\leq k \leq n/2} q^{u(k)} (1-q^{-\ffr{r-1}{k}}) & <
    q^{u(1)} (1-q^{-r}) + \sum_{2\leq k
      \leq n/2}
    q^{u(k)}  \label{eq:24} \\
    & < q^{u(1)}  (1-q^{-r}) + \frac{q^{u(2)}}{1-q^{-1}}  \\
    & = q^{u(1)}(1-q^{-r}) \Bigl( 1 +
      \frac{q^{-u(1)+u(2)}}{(1-q^{-1})(1-q^{-r})}
    \Bigr).
  \end{align}
  Since $u(1)=b_{r,n-1}+r-1$ and
  $-u(1)+u(2)=-\ffr{r-1}{n-1}+\ffr{r-1}{2}$, we conclude that
  \begin{align}
    \# \Rone{r,n} & \leq
    \frac{q^{b_{r,n-1}+r-1}(1-q^{-r})}{(1-q^{-1})^2} \Bigl( 1 +
      \frac{q^{-\ffr{r-1}{n-1}+\ffr{r-1}{2}}}{(1-q^{-1})(1-q^{-r})}
    \Bigr) \\
    & = \rho_{r,n}(q) \Bigl( 1 +
      \frac{q^{-\ffr{r-1}{n-1}+\ffr{r-1}{2}}}{(1-q^{-1})(1-q^{-r})}
    \Bigr) \\
    & < \rho_{r,n}(q) ( 1 + 3q^{-\ffr{r-1}{n-1}+\ffr{r-1}{2}} ). \label{eq:36}
  \end{align}
This proves \eqref{eq:33} and we proceed with \eqref{eq:34}.  Using
\eqref{eq:36}, we have
  \begin{align}
    \#\Ione{r,n} & = \# \Pone{r,n} - \# \Rone{r,n} \\
    & \geq \# \Pone{r,n} \Bigl(1 - \rho_{r,n}(q)\frac{1+ 3 q^{-b_{r-1,n-1}+b_{r-1,2}}}{\# \Pone{r,n}}\Bigr) \\
    & = \# \Pone{r,n} \Bigl( 1 -  q^{-b_{r-1,n}+r} \frac{(1+3
        q^{-b_{r-1,n-1}+b_{r-1,2}})(1-q^{-r})}{(1-q^{-1})(1-q^{-b_{r-1,n}})}
    \Bigr).
  \end{align}
We observe that the exponent
  $-b_{r-1,n-1}+b_{r-1,2}$ is decreasing in $r$ and $n$ for $n \geq 4$. It is
  furthermore always negative and hence the fraction $(1+3 q^{-b_{r-1,n-1}+b_{r-1,2}})/(1-q^{-b_{r-1,n}})$ is also
  decreasing in $q$. Therefore it achieves its maximal value for $n=4$, $r=2$ and $q=2$, yielding $80/31<3$ as upper bound and proving \eqref{eq:34}.
  For the last argument, we need
  \eqref{eq:34} also for $n=3$;  this follows from
  \autoref{thm:R-from-gen}.

  We conclude with the proof of \eqref{eq:35}.  The subset $\{g \cdot h \colon g \in \Pone{r,1}, h \in \Ione{r,n-1} \} \subset \Rone{r,n,k}$ has size
  $\# \Pone{r,1} \cdot \# \Ione{r,n-1}$.  With \eqref{eq:34}, we find
  \begin{align}
    \# \Rone{r,n} & \geq \# \Pone{r,1} \cdot \# \Ione{r,n-1}  \\
    & \geq q^{b_{r,1}-1}\frac{1-q^{-r}}{1-q^{-1}} \cdot
    \# \Pone{r,n-1} \Bigl(1- 3 q^{-b_{r-1,n-1}+r} \frac{1-q^{-r}}{1-q^{-1}}\Bigr)  \\
 & = \rho_{r,n}(q) (1-q^{-b_{r-1,n-1}})  \Bigl(1- 3 q^{-b_{r-1,n-1}+r} \frac{1-q^{-r}}{1-q^{-1}}\Bigr)  \\
    & \geq \rho_{r,n}(q)  \Bigl(1- 3 q^{-b_{r-1,n-1}+r} \frac{1-q^{-r-1}}{1-q^{-1}}\Bigr).
  \end{align}
  We combine the upper and lower bounds \eqref{eq:33}
  and \eqref{eq:35}.  The maximum
  of the bounds on the relative error term is
  \[
  \max\Bigl(3 q^{-r(r-1)/2} (1-q^{-r-1}),\frac{1}{1-q^{-r}} \Bigr) \cdot
  \frac{q^{-\ffr{r-1}{n-1}+\ffr{r-1}{2}}}{1-q^{-1}} =  \frac{q^{-\binom{r+n-2}{r-1}+r(r+1)/2}}{(1-q^{-1})(1-q^{-r})}
  \]
and the observation $(1-q^{-1})(1-q^{-r}) \leq 8/3$ concludes the proof.
\end{proof}

The approach of this section also works, with
minor modifications, for $n<4$ and can provide a stand-alone
proof of \autoref{thm:red}, without recourse to \autoref{thm:R-from-gen}.

\begin{figure}
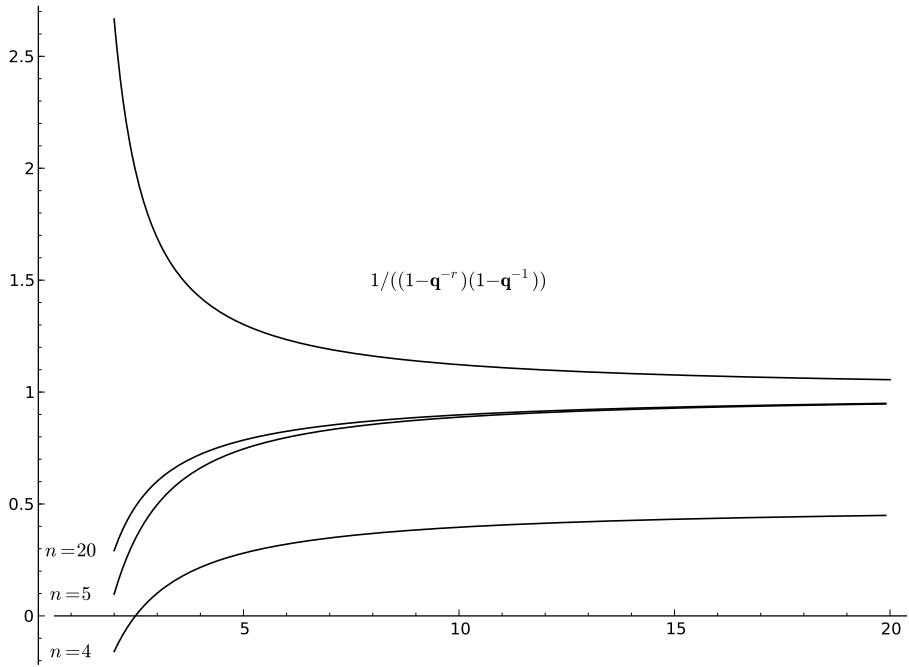

\centering
\sageplot[width=\textwidth]{RR_vs_rho_plot(2,[4,5,20])}
\caption{The normalized relative error in \autoref{thm:R-from-gen} for
  $r=2$.}
\label{fig:RR_vs_rho}
\end{figure}

\autoref{fig:RR_vs_rho} shows plots of $(\scoeff{R_n}(\qvar)-\rho_{r,n}(\qvar))/(\rho_{r,n}(\qvar)
  \qvar^{-\binom{r+n-2}{r-1}+r(r+1)/2})$ for $r=2$ and $n=4,5,20$ as we substitute for
  $\qvar$ real numbers from $2$ to $20$.  \autoref{thm:red} says that
  the values are absolutely at most $1/((1-\qvar^{-r})(1-\qvar^{-1}))$.  \autoref{thm:R-from-gen} indicates a bound of $1/2 + o(1)$ for $n=4$ and $1+o(1)$ for $n>4$, but without explicit error estimate.

According to \eqref{eq:16}, the bound on the absolute value of the relative error for $n \geq 4$ is
\begin{equation}
  \label{eq:7}
\frac{q^{-b_{r-1,n-1}+b_{r-1,2}}}{(1-q^{-1})(1-q^{-r})}.
\end{equation}
For $n>4$, this is at most $2/3$.  For $n=4$, we can drop the factor $1-q^{-1}$, since the sum in \eqref{eq:24} consists only of a single summand and the estimate by a geometric sum is not necessary.  This shows that also for $n=4$, the relative error is at most $2/3$.

\begin{remark}
  \label{rem:exp_decay}
  How close is our relative error estimate to being exponentially
  decaying in the input size? The usual dense representation of a
  polynomial in $r$ variables and of degree $n$ requires $b_{r,n} =
  \binom{r+n}{r}$ monomials, each of them equipped with a coefficient
  from $\FF_q$, using about $\log_2 q$ bits. Thus the total input size
  is about $\log_2 q \cdot b_{r,n}$ bits.  This differs from $\log_2 q
  \cdot (b_{r-1,n-1}-b_{r-1,2})$ by a factor of
  \[
  \frac{b_{r,n}}{b_{r-1,n-1}-b_{r-1,2}} <
  \frac{b_{r,n}}{\frac{1}{2}b_{r-1,n-1}} = \frac{2(n+r)(n+r-1)}{nr}.
  \]
  Up to this polynomial difference (in the exponent), the relative
  error is exponentially decaying in the bit size of the input, that
  is, $(\log q)$ times the number of coefficients in the usual dense
  representation.  In particular, it
  is exponentially decaying in any of the parameters $r$, $n$, and $\log_2
  q$, when the other two are fixed.
\end{remark}

These bounds fit well into the picture described in Section 2 of
\cite{gat08-incl-gat07} for $r=2$.  The family of functions described
there approximates the quotient $\# R_{2,n} / \# P_{2,n}$ (using our
notation).  If we compare them to $\rho_{r,2}(q) / \# \Pone{2,n}$ we
find that they differ only by the factor $1-q^{-n-1}$, which tends to
$1$ as $n$ and $q$ increase.  Our bound  $3 q^{-n +3}$ on the relative error
for $r=2$ and $n \geq 4$ is only slightly larger than the bound
$2q^{-n+3}$ in Theorem~2.1(ii) of the paper cited.

The following provides some handy bounds.

\begin{corollary}
\label{cor:handy_R}
For $r, q \geq 2$, and $n\geq 5$, we have
  \begin{gather}
\frac{1}{4} q^{\binom{r+n-1}{r}+r-1} \leq \# \Rone{r,n}(\FF_{q}) \leq 6 
 q^{\binom{r+n-1}{r}+r-1},     \label{eq:17} \\
\frac{1}{4} q^{-\binom{r+n-1}{r-1}+r} \leq \frac{\# \Rone{r,n}(\FF_{q})}{\# \Pone{r,n}(\FF_{q})} \leq 3 
q^{-\binom{r+n-1}{r-1}+r}.
  \end{gather}
\end{corollary}


We conclude this section with bounds for the number of irreducible polynomials.

\begin{corollary}
  \label{cor:irred}
  Let $r, q\geq 2$, and $\rho_{r,n}$ as in \autoref{thm:R-from-gen}.  We have
  \begin{equation}
\label{eq:102}
\# \Pone{r,n} (\FF_q) - 2 \rho_{r,n}(q) \leq  \# \Ione{r,n} (\FF_q) \leq \# \Pone{r,n} (\FF_q) ,
  \end{equation}
  and more precisely
  \begin{align*}
    \# \Ione{r,1} (\FF_{q}) & = \# \Pone{r,1} (\FF_{q}), \\
    \# \Ione{r,2} (\FF_q) & = \# \Pone{r,2} (\FF_q) - \frac{
      \rho_{r,2}(q)}{2} \cdot
    (1-q^{-r-1}), \\
    \left| \# \Ione{r,3} (\FF_q)- (\# \Pone{r,3}(\FF_q) -
        \rho_{r,3}(q)) \right| & \leq \rho_{r,3}(q) \cdot
    q^{-(r-1)r/2}, \\
    \intertext{and for $n \geq 4$} \left| \# \Ione{r,n} (\FF_q)-
    (\# \Pone{r,n}(\FF_q) - \rho_{r,n}(q)) \right| & \leq
    \rho_{r,n}(q) \cdot 3 q^{-\binom{r+n-2}{r-1}+r(r+1)/2}.
  \end{align*}
\end{corollary}

\begin{proof}
  The more precise statements follow directly from \autoref{thm:red} by
  application of $\# \Pone{r,n}(\FF_q) = \# \Rone{r,n}(\FF_q) + \#
  \Ione{r,n}(\FF_q)$.  These imply the first claim for $n<4$.  For $n
  \geq 4$, the relative error in \eqref{eq:16} is at most $2/3<1$ as
  remarked after the proof of \autoref{thm:red} and this concludes the
  proof of \eqref{eq:102}.
\end{proof}

\section{Powerful polynomials}
\label{sec:powerful}

For an integer $s \geq 2$, a polynomial is called
\emph{$s$-powerful} if it is divisible by the $s$th power of some
nonconstant polynomial, and \emph{$s$-powerfree} otherwise; it is
\emph{squarefree} if $s=2$. Let
\begin{align*}
  \Qone{r,n,s}   (\F) & = \{f\in \Pone{r,n}  (\F) \colon f \text{ is
    $s$-powerful}\}, \\
  \Sone{r,n,s} (\F) & = \Pone{r,n}(\F) \setminus \Qone{r,n,s}(\F).
\end{align*}
As in the previous section, we restrict our attention to a finite field
$\F = \FF_q$, which we omit from the notation.

For the approach by generating functions, we consider the
combinatorial classes $\ccQ = \bigcup_{n \geq 0} \Qone{r,n,s}$ and
$\ccS = \ccP \setminus \ccQ$, where the explicit reference to
$r$ and $s$ is omitted.  Any monic polynomial $f$ factors uniquely as $f=g\cdot
h^s$ where $g$ is a monic $s$-powerfree polynomial and $h$ an arbitrary
monic polynomial, hence
\begin{equation}
  \label{eq:61}
\nP = \nS \cdot \nP ( z^{s})
\end{equation}
and  by definition $\nQ = \nP
- \nS$ for the generating functions of $\ccS$ and $\ccQ$,
respectively.  For univariate polynomials, \cite{car32} derives \eqref{eq:61}
directly from generating functions to prove the counting formula which
we reproduce in \eqref{eq:52}.  \citet[Section~1.1]{flagou01} use
\eqref{eq:61} for $s=2$ to count univariate squarefree polynomials,
see also
\citet[Note~I.66]{flased09}.  A corresponding Maple program to compute the coefficients of $\nQ$ is shown in \autoref{fig:algopower}. It was used to compute
$\# \Qone{2,n,2} (\FF_{q})$ for $n\leq 6$ in \citet[Table~3.1]{gat08-incl-gat07}.  We extend this in \autoref{tablaQ}.

\begin{figure}
\centering
\begin{Verbatim}[frame=single]
spowerfreesGF:=proc(z,N,r,s) local i: option remember:
    convert(taylor(allpolysGF(z,N,r)/allpolysGF(z^s,N,r),
        z,N+1),polynom):
end:

spowerfulsGF:=proc(z,N,r,s) option remember:
    allpolysGF(z,N,r)-spowerfreesGF(z,N,r,s):
end:

spowerfuls:=proc(n,r,s)
    coeff(sort(expand(spowerfulsGF(z,n,r,s))),z^n):
end:
\end{Verbatim}
\caption{Maple program to compute the number of monic
$s$-powerful polynomials in $r$ variables of degree $n$.}
\label{fig:algopower}
\end{figure}

\begin{table}\footnotesize
\centering
\sagestr{Q_table(2,9,3)}

\sagestr{Q_table(3,6,2)}

\sagestr{Q_table(3,7,3)}
  \caption{Exact values of $\# \Qone{r,n,s} (\FF_{q})$ for small values
of $r,n, s$.}
\label{tablaQ}
\end{table}

As in \autoref{pro:R_exact_by_recursion}, this approach quickly leads to explicit formulas.
\begin{theorem}
\label{pro:Q_exact_by_recursion}
For $r \geq 1$, $q,s \geq 2$, $\nP_{n}$ as in \eqref{eq:44}, and $M_{n}$ as in \eqref{eq:10}, we have
  \begin{align}
\nS_{n} & = \sum_{\substack{0 \leq i \leq n/s \\ j \in M_{i}}} (-1)^{\abs{j}}  \nP_{j_{1}} \nP_{j_{2}} \cdots \nP_{j_{\abs{j}}} \nP_{n-is}, \\
    \nQ_{n} & = - \sum_{\substack{1 \leq i \leq n/s \\ j \in M_{i}}} (-1)^{\abs{j}}  \nP_{j_{1}} \nP_{j_{2}} \cdots \nP_{j_{\abs{j}}} \nP_{n-is} \label{eq:71}.
  \end{align}
\end{theorem}

\begin{proof}
We consider the original power series $\ngenfun{F} = 1 - \nP = - \sum_{i \geq 1} \nP_{i} z^{i}$ and express \eqref{eq:61} as
\begin{align}
  \nS & = \nP \cdot \sum_{i \geq 0} \ngenfun{F} (z^{s})^{i} \\
& = \sum_{k \geq 0} \nP_{k} z^{k}  \cdot \sum_{i \geq 0} \Bigl(- \sum_{j \geq 1} \nP_{j} z^{js} \Bigr)^{i}.
\end{align}
Comparison of coefficients provides us with
\begin{align}
\nS_{n} & = \sum_{\substack{0 \leq i \leq n/s \\ j \in M_{i}}} (-1)^{\abs{j}}  \nP_{j_{1}} \nP_{j_{2}} \cdots \nP_{j_{\abs{j}}} \nP_{n-is},
\end{align}
and the claim for $\nQ_{n} = \nP_{n} - \nS_{n}$ follows.
\end{proof}

For $r=1$, we have $\nP_j = q^j$ and for any composition $j_{1} + j_{2} +  \dots + j_{k}$ of $i$ in \eqref{eq:71}
\begin{equation}
  \nP_{j_{1}} \nP_{j_{2}} \cdots \nP_{j_{k}} \nP_{n-is} = q^{n-(s-1)i}.
\end{equation}
Moreover, since
\begin{equation}
\sum_{k\geq 1} (-1)^k \binom{i-1}{k-1} = - \binom{0}{i-1} = \begin{cases}
-1 & \text{if $i=1$}, \\
0 & \text{if $i\geq2$},
\end{cases}
\end{equation}
see \citet[p.~167]{graknu89}, we have in the univariate case
\begin{equation}
\label{eq:52}
\nQ_{n}
= - \sum_{\substack{1 \leq i \leq n/s \\
k \geq 1 }} (-1)^k \binom{i-1}{k-1} q^{n-(s-1)i} =
\begin{cases}
0  & \text{if $n < s$}, \\
q^{n-s+1} & \text{if $n \geq s$},
\end{cases}
\end{equation}
as shown by \citet[Section 6]{car32}.

\begin{figure}
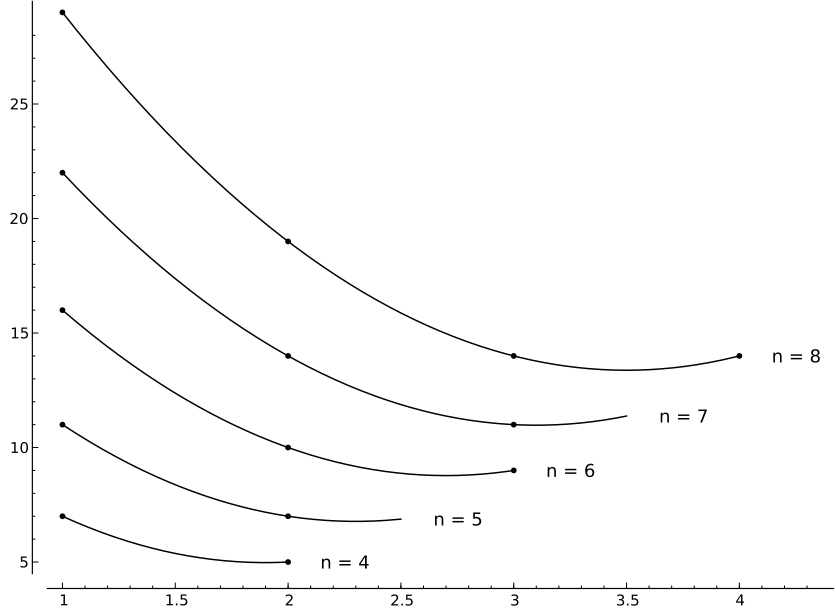

\centering
\sageplot[width=\textwidth]{v_plot(2,srange(4,9),2)}
\caption{Graphs of $v_{2,n,2}(k)$ on $[1,n/2]$ as $n$ runs from $4$ to
  $8$.  The dots represent the values at integer arguments.}
\label{fig:v}
\end{figure}

To study the asymptotic behavior of $\nS_{n}$ and $\nQ_{n}$ for $r \geq 2$ we again deviate from the standard approach and move to power series in $\Qqz$.  With $\sP$ from \eqref{eq:5}, we define $\sS, \sQ \in \Qqz$ by
  \begin{align}
    \sP & = \sS \cdot \sP(z^{s}), \label{def:S} \\
    \sQ & = \sP - \sS.
\end{align}
This is well-defined, since $\sP (z^{s})$ has constant term 1 and is
  therefore invertible.
By construction, we have
\begin{equation}\label{eq:25}
  \begin{split}
    \scoeff{S}_n (q) & = \# \Sone{r,n,s} (\FF_q), \\
    \scoeff{Q}_n (q) & = \# \Qone{r,n,s} (\FF_q).
  \end{split}
\end{equation}
To study the asymptotic behavior, we examine $\scP_k \cdot
  \scP_{n-sk}$. Let
\begin{align}
  v_{r,n,s} (k) & = \degq (\scP_k \cdot
  \scP_{n-sk}) \label{eq:2010} \\
  & = (r+k)\ttf{r} / r!  + (r+n-sk)\ttf{r}/ r! -2
\end{align}
and consider $v_{r,n,s}(k)$ as a function of a real variable $k$ (\autoref{fig:v}).  In contrast to $u(k)$ from Section~\ref{sec:gen}, this function is not monotone in $k$.

\begin{lemma}
  \label{lem:conv}
  Let $r,n,s,q \geq 2$.
  \begin{ronumerate}
  \item \label{item:3} The function $v_{r,n,s}(k)$ is convex for $1\leq k \leq n/s$.
\item \label{lem:conv:ii} For all integers $k$ with $2 \leq k \leq
  n/s$, we have
\[
v_{r,n,s} (1) > v_{r,n,s}(k).
\]
\item \label{lem:62} For all integers $k$ with $3 \leq k \leq n/s$, we have
\begin{align}
    v_{r,n,s}(2) & > v_{r,n,s}(k) && \text{if $(n,s) \neq (6,2)$}. \\
\intertext{Furthermore,}
    v_{r,6,2}(2) & < v_{r,6,2}(3) && \text{if $r \geq 3$}, \\
    v_{2,6,2}(2) & = v_{2,6,2}(3) +1. &&
\end{align}
\item  \label{cor:geosum2}
If $(n,s)\neq (6,2)$, then
  \begin{equation}
    \label{eq:geosum2}
    \sum_{2\leq k \leq n/s} q^{v_{r,n,s}(k)} \leq
    \frac{2q^{v_{r,n,s}(2)}}{1-q^{-1}}.
  \end{equation}
  \end{ronumerate}
\end{lemma}

\begin{proof}
  We switch to the affine transformation
  \begin{align}
    \bar{v}(k) & =  r! \cdot \Bigl(v_{r,n,s}(k)+2\Bigr) \\
    & = (r+k)\ttf{r} + (r+n-sk)\ttf{r},
  \end{align}
  which exhibits the same behavior as $v_{r,n,s}$ concerning convexity and maximality.
  \begin{ronumerate}
  \item We have
  \begin{equation}
    \bar{v}''(k) =\sum_{\substack{1 \leq i,j \leq r \\ i \neq j}} \Bigl(
      \frac{(r+k)\ttf{r}}{(i+k)(j+k)} +
      \frac{s^{2}(r+n-sk)\ttf{r}}{(i+n-sk)(j+n-sk)} \Bigr) >0.
  \end{equation}
\item For $n < 2s$, there is nothing to prove.  For $n \geq 2s$, we find $n \geq s+2 \geq s+1+1/(s-1)$ and for all $i$
\begin{align}
&   (i+ n-s) - (i+ n/s) \geq 0, \\
& (r+n-s)\ttf{r} - (r+n/s)\ttf{r} \geq 0, \\
\bar{v}(1)-\bar{v}(n/s) & = (r+1)! + (r+n-s)\ttf{r} -(r+n/s)\ttf{r} -r! \\
     & =(r+n-s)\ttf{r} - (r+n/s)\ttf{r} +  r\cdot r!> 0.
   \end{align}
With the convexity of $\bar{v}$, this suffices.

\item Analogously to \ref{lem:conv:ii}, it is
sufficient to prove $\bar{v}(2) > \bar{v}( n/s )$ for
$(n,s)\neq(6,2)$. If $n\geq 2s^2/(s-1)$, then $n-2s \geq n/s$, so that for all $i$
\[
(i + n-2s) - (i+n/s) \geq 0
\]
and hence
  \begin{align*}
    \bar{v}(2) - \bar{v}(n/s) &= (r+2)!/2 + (r+n-2s)\ttf{r} - (r+n/s)\ttf{r} - r!  \\
    & > (r+n-2s)\ttf{r} - (r+n/s)\ttf{r} \geq 0.
  \end{align*}
  If $n<2s^2/(s-1)$, then $n/s < 3$ for $s \geq 3$ or $n < 6$ and
  there is nothing to prove.  Finally, the three conditions $n<2s^2/(s-1)$, $s=2$, and $n
  \geq 6$ enforce $6\leq n < 8$, and we compute directly
  \[
  v_{r,7,2}(2)-v_{r,7,2}(3) = \frac{1}{2} r (r+1) >0,
  \]
  \[
  v_{r,6,2}(2)-v_{r,6,2}(3)=-\frac{1}{6}(r-3)(r+1)(r+2)-1
  \begin{cases}
    =1 & \text{ if } r=2, \\
    <0 & \text{ if } r\geq 3.
  \end{cases}
  \]

\item The maximal value of the integer sequence $ v_{r,n,s}(k)$ for $2 \leq k \leq n/s$ is $v_{r,n,s}(2)$ by \ref{lem:62}.  Each value is taken at most twice, due to \ref{item:3}, and we can bound the sum by twice a geometric sum as
\begin{equation}
    \sum_{2\leq k \leq n/s} q^{v_{r,n,s}(k)} \leq 2 q^{v_{r,n,s}(2)}  \sum_{k\geq 0} q^{-k} =
    \frac{2q^{v_{r,n,s}(2)}}{1-q^{-1}}.  \mbox{\qedhere}
\end{equation}
\end{ronumerate}
\end{proof}
The approach by generating functions now yields the following result.  Its ``general'' case is \ref{item:13}.  We give exact expressions in special cases, namely for $ n < 3s$ in \ref{item:17} and for $(n,s)=(6,2)$ in \ref{item:18}, which also apply when we substitute the size $q$ of a finite field $\FF_{q}$ for $\qvar$.

\begin{theorem}
\label{thm:Q_by_gen}
  Let $r, s \geq 2$, $n\geq 0$,  and
  \begin{align}
  \eta_{r,n,s} (\qvar) & = \qvar^{\binom{r+n-s}{r} + r-1}
  \frac{(1-\qvar^{-r}) (1-\qvar^{-\binom{r+n-s-1}{r-1}})}{(1-\qvar^{-1})^2} \in
  \QQ (\qvar), \label{eq:89} \\
  \delta & = \binom{r+n-s}{r}-\binom{r+n-2s}{r} - \frac{r(r+1)}{2}.
  \end{align}
  \begin{ronumerate}
  \item \label{item:4} If $n \geq 2s$, then $\delta \geq r$.
  \item \label{item:17} \begin{equation}
      \label{eq:50}
      \scoeff{Q}_{n} = \begin{cases}
0 & \text{for $n<s$,}\\
\eta_{r,n,s}(\qvar) & \text{for $s \leq n<2s$,} \\
\eta_{r,n,s} (\qvar) \biggl( 1 + \qvar^{-\delta} \cdot \frac{1-\qvar^{-\binom{n+r-2s-1}{r-1}}}{1-\qvar^{-\binom{n+r-s-1}{r-1}}} & \\
\quad \cdot \Bigl( \frac{1-\qvar^{-r(r+1)/2}}{1-\qvar^{-r}} - \qvar^{-r(r-1)/2}\frac{1-\qvar^{-r}}{1-\qvar^{-1}}\Bigr)\biggr) & \text{for $2s \leq n<3s$.}
\end{cases}
\end{equation}

  \item \label{item:18} For $(n,s)=(6,2)$ and $r \geq 2$, we have
    \begin{align}
      \scoeff{Q}_{6}  & = \eta_{r,6,2}(\qvar) \bigg(1 + \qvar^{- \binom{r+3}{4}-r+1} \cdot \Big( \qvar^{-1}  \tfrac{(1-\qvar^{-1})\bigl(1-\qvar^{-\binom{r+2}{3}}\bigr)}{(1-\qvar^{-r})\bigl(1-\qvar^{-\binom{r+3}{4}}\bigr)} \\
& \quad +\qvar^{-(r^{3}-7r+6)/6} \tfrac{(1-\qvar^{-r(r+1)/2})^{2}}{(1-\qvar^{-r})\bigl(1-\qvar^{-\binom{r+3}{4}}\bigr)} \\
& \quad - \qvar^{-(r^{3}+3r^{2}-10r+6)/6} \tfrac{(1-\qvar^{-r})(1-\qvar^{-r(r+1)/2})}{(1-\qvar^{-1})\bigl(1-\qvar^{-\binom{r+3}{4}}\bigr)} \\
& \quad  -2\qvar^{-(r^{3}+3r^{2}+4r-6)/6} \tfrac{1-\qvar^{-r(r+1)/2}}{1-\qvar^{-\binom{r+3}{4}}} \\
& \quad + \qvar^{-(r^{3}+6r^{2}-7r+6)/6} \tfrac{(1-\qvar^{-r})^{2}}{(1-\qvar^{-1})\big(1-\qvar^{-\binom{r+3}{4}}\big)}\Big)\bigg) \\
& = \eta_{r,6,2} (\qvar) \bigl( 1 + \qvar^{-\delta+(r-2)(r-1)(r+3)/6}(1+O(\qvar^{-1}))\bigr). \label{eq:86}
    \end{align}
  \item\label{item:13} For $n \geq 2s$ and $(n,s)\neq (6,2)$, we have
    \begin{equation}
\label{eq:49}
      \scoeff{Q}_{n} =   \eta_{r,n,s}(\qvar) \big(1 + \qvar^{-\delta} (1+    O(\qvar^{-1}))\big).
    \end{equation}
  \end{ronumerate}
\end{theorem}
\begin{proof}
\begin{ronumerate}
\item If $n \geq 2s$, then
  \begin{equation}
    \delta  \geq \binom{r+s}{r}-1-\frac{r(r+1)}{2}
 \geq \binom{r+2}{r}-1-\frac{r(r+1)}{2} = r.
  \end{equation}
  \item The exact formulas of \autoref{pro:Q_exact_by_recursion} yield
\begin{align}
\scoeff{Q_{n}} & = 0 && \text{for $n<s$}, \\
  \scoeff{Q_n} & = \scP_1 \scP_{n-s} = \eta_{r,n,s}(\qvar) && \text{for $s \leq n < 2s$},
\end{align}
and for $2s \leq n < 3s$,
\begin{align}
 \scS_{n}  & = \scP_{n} - \scP_{1}\scP_{n-s} - (\scP_{2}-\scP_{1}^{2})\scP_{n-2s}, \label{eq:46} \\
\scQ_{n} & = \scP_{1}\scP_{n-s} + (\scP_{2} - \scP_{1}^{2})\scP_{n-2s} \\
& = \eta_{r,n,s} (\qvar) \Bigl( 1 + \frac{\scP_{2}\scP_{n-2s}}{\scP_{1}\scP_{n-s}}\Bigl(1-\frac{\scP_{1}^{2}}{\scP_{2}}\Bigr)\Bigr) \\
& = \eta_{r,n,s} (\qvar) \bigg( 1 + \qvar^{-\delta} \frac{1-\qvar^{-\binom{n+r-2s-1}{r-1}}}{1-\qvar^{-\binom{n+r-s-1}{r-1}}}  \\
& \quad \quad \cdot \Big( \frac{1-\qvar^{-r(r+1)/2}}{1-\qvar^{-r}} - \qvar^{-r(r-1)/2}\frac{1-\qvar^{-r}}{1-\qvar^{-1}}\Big)\bigg),
\end{align}
where  $\delta = - \degq (\scP_{2}\scP_{n-2s}/(\scP_{1}\scP_{n-s}))$.

\item
For $s=2$, we evaluate \eqref{eq:71} for
\begin{align}
\scQ_{6} & =  \scP_{1}\scP_{4} + \scP_{3} + \scP_{2}^{2} - \scP_{1}^{2}\scP_{2} - 2 \scP_{1}\scP_{2} + \scP_{1}^{3} \\
& = \eta_{r,6,2} (\qvar) (1 + ( \scP_{3} + \scP_{2}^{2} - \scP_{1}^{2}\scP_{2} - 2\scP_{1}\scP_{2} + \scP_{1}^{3})/(\scP_{1}\scP_{4})) \\
& = \eta_{r,6,2}(\qvar) \biggl(1 + \qvar^{-v_{r,6,2}(1)+v_{r,6,2}(3)+1}
  \\
& \quad \cdot \Bigl( \qvar^{-1}
  \tfrac{(1-\qvar^{-1})(1-\qvar^{-b_{r-1,3}})}{(1-\qvar^{-r})(1-\qvar^{-b_{r-1,4}})}
 \label{eq:109} \\
& \quad +\qvar^{-(r^{3}-7r+6)/6} \tfrac{(1-\qvar^{-r(r+1)/2})^{2}}{(1-\qvar^{-r})(1-\qvar^{-b_{r-1,4}})} \\
& \quad - \qvar^{-(r^{3}+3r^{2}-10r+6)/6} \tfrac{(1-\qvar^{-r})(1-\qvar^{-r(r+1)/2})}{(1-\qvar^{-1})(1-\qvar^{-b_{r-1,4}})}   \\
& \quad -2\qvar^{-(r^{3}+3r^{2}+4r-6)/6} \tfrac{1-\qvar^{-r(r+1)/2}}{1-\qvar^{-b_{r-1,4}}} \\
& \quad + \qvar^{-(r^{3}+6r^{2}-7r+6)/6}
    \tfrac{(1-\qvar^{-r})^{2}}{(1-\qvar^{-1})(1-\qvar^{-b_{r-1,4}})}\Bigr)\biggr) \label{eq:51}
\\
& = \eta_{r,6,2} (\qvar) \big( 1 + \qvar^{-\delta+(r-2)(r-1)(r+3)/6}(1+O(\qvar^{-1}))\big), \label{eq:54}
\end{align}
since the sum \eqref{eq:109}--\eqref{eq:51} has nonpositive degree in $\qvar$
and $-v_{r,6,2}(1) +v_{r,6,2}(3)+1=
 - \binom{r+3}{4}-r+1 =
-\delta + (r-2)(r-1)(r+3)/6$.

\item Finally, for $n \geq 2s$ and $(n,s)\neq (6,2)$, we claim
\begin{equation}
\label{eq:47}
  \scS_{n} = \scP_{n} - \scP_{1}\scP_{n-s} - \scP_{2}\scP_{n-2s}(1+ O (\qvar^{-1})).
\end{equation}
This implies immediately
\begin{equation}
\label{eq:55}
  \begin{split}
    \scS_{n} & = \scP_{n} - \scP_{1}\scP_{n-s} (1+ O (\qvar^{-1})) \\
    & =\scP_{n}(1+ O (\qvar^{-1})),
  \end{split}
\end{equation}
by Lemmas~\ref{lem:conv}~\ref{lem:conv:ii} and \ref{lem:deg_q}~\ref{item:degi}, respectively. We already have \eqref{eq:47} for $ 2s \leq n < 3s$ from \eqref{eq:46} by \autoref{lem:deg_q}~\ref{item:degi}.  We also have \eqref{eq:55} for $(n,s)=(6,2)$ from \eqref{eq:54}.  This is enough to obtain inductively
\begin{align}
  \scS_{n} & = \scP_{n} - \sum_{1 \leq i \leq n/s} \scS_{n-is} \scP_{i} \\
  &= \scP_{n} - \scP_{1}\scS_{n-s} - \sum_{2 \leq i \leq n/s} \scP_{i}\scS_{n-is} \\
  & = \scP_{n}-\scP_{1}\bigl(\scP_{n-s} - \scP_{1}\scP_{n-2s}(1+O(\qvar^{-1}))\bigr) - \sum_{2  \leq i \leq n/s} \scP_{i}\scP_{n-is}( 1 + O( \qvar^{-1})) \\
  & = \scP_{n} - \scP_{1}\scP_{n-s} +
  \scP_{1}^{2}\scP_{n-2s}(1+O(\qvar^{-1}))
  -\scP_{2}\scP_{n-2s}(1+ O(\qvar^{-1})) \\
& = \scP_{n} - \scP_{1}\scP_{n-s} - \scP_{2}\scP_{n-2s}(1+ O (\qvar^{-1})),
\end{align}
using \autoref{lem:conv}~\ref{lem:62} for $(n,s)\neq(6,2)$ and
\autoref{lem:deg_q}~\ref{item:degi}.  We conclude with
\begin{align}
  \scQ_{n} & = \scP_{1}\scP_{n-s} + \scP_{2}\scP_{n-2s}(1+ O (\qvar^{-1})) \\
& =  \eta_{r,n,s} (\qvar) ( 1 + \qvar^{-\delta} ( 1 + O( \qvar^{-1})) )
\end{align}
by the definition of $\eta_{r,n,s}(\qvar) = \scP_{1}\scP_{n-s}$ and
$\delta = - \degq (\scP_{2}\scP_{n-2s}/(\scP_{1}\scP_{n-s}))$,
respectively. \qedhere
\end{ronumerate}
\end{proof}
For $r\geq 3$, we can replace $1+O(\qvar^{-1})$ in \eqref{eq:86} by $\qvar^{-1}+O(\qvar^{-2})$.

In the following, the combinatorial approach replaces the asymptotic $1+O(\qvar^{-1})$ of \eqref{eq:49} with an explicit bound of $6$ in \eqref{eq:45}.  We consider for integers $1\leq k \leq n/s$ the sets
\begin{equation}
  \Qone{r,n,s,k}(F) = \{ g \cdot h^{s} \colon g \in \Pone{r,n-sk}, h \in \Pone{r,k}\} \in \Pone{r,n} (F)
\end{equation}
and have
\begin{equation}
  \label{eq:prop_sigma}
  \Qone{r,n,s} (F)  = \bigcup_{1\leq k \leq n/s} \Qone{r,n,s,k} (F).
\end{equation}
For $n<3s$ the exact formula \eqref{eq:50} of \autoref{thm:Q_by_gen}~\ref{item:17} applies.  We provide explicit bounds for $n \geq 3s$.

\begin{theorem}
\label{thm:Q_by_map}
Let $r,s,q\geq 2$, $n \geq 0$, and
  \begin{align}
  \eta_{r,n,s} (\qvar) & = \qvar^{\binom{r+n-s}{r} + r-1}
  \frac{(1-\qvar^{-r}) \bigl(1-\qvar^{-\binom{r+n-s-1}{r-1}}\bigr)}{(1-\qvar^{-1})^2} \in
  \QQ (\qvar), \\
  \delta & = \binom{r+n-s}{r}-\binom{r+n-2s}{r} - \frac{r(r+1)}{2}
  \end{align}
as in \autoref{thm:Q_by_gen}.
\begin{ronumerate}
\item \label{item:8} For $(n,s) = (6,2)$, we have $\delta = r(r+1)(r^{2}+9r+2)/24$ and
\begin{equation}
\label{eq:53}
  \left| \# \Qone{r,6,2} (\FF_q) - \eta_{r,6,2}(q) \right| \leq
  \eta_{r,6,2}(q) \cdot 2 q^{-\delta+(r-2)(r-1)(r+3)/6}.
\end{equation}
  \item \label{item:9} For $n \geq 3s$ and $(n,s) \neq (6,2)$, we have
\begin{equation}
\label{eq:45}
  \left| \# \Qone{r,n,s} (\FF_q) - \eta_{r,n,s}(q) \right| \leq
  \eta_{r,n,s}(q) \cdot 6 q^{-\delta}.
\end{equation}

\end{ronumerate}
\end{theorem}

\begin{proof}
We omit the argument $\FF_{q}$ from the notation. Considering only the positive and negative summands  of \eqref{eq:51}, respectively, we find
\begin{align}
  \# \Qone{r,6,2} & \leq \eta_{r,6,2} (q) (1+ 2q^{-\delta+(r-2)(r-1)(r+3)/6}), \label{eq:9} \\
  \# \Qone{r,6,2} & \geq \eta_{r,6,2} (q) (1 - q^{-\delta+(r-2)(r-1)(r+3)/6}),
\end{align}
which proves \ref{item:8}.

For the general case \ref{item:9}, we claim
  \begin{align}
    \# \Qone{r,n,s} & \leq \eta_{r,n,s}(q) \Bigl( 1 + \frac{16}{3}q^{-\delta}\Bigr) &&\text{for $(n,s)\neq(6,2)$}, \label{eq:31} \\
    \# \Qone{r,n,s} & \geq \eta_{r,n,s}(q) \Bigl( 1 - \frac{7}{2}q^{-\delta-r(r-1)/2}\Bigr) &&\text{for $n \geq 3s$}. \label{eq:32}
  \end{align}

For \eqref{eq:31}, we find from \eqref{eq:prop_sigma}
\begin{align}
  \# \Qone{r,n,s} & \leq \sum_{1\leq k \leq n/s} \# \Qone{r,n,s,k} \leq \sum_{1\leq k \leq n/s} \# \Pone{r,n-sk} \cdot \# \Pone{r,k}
  \\
  & = \sum_{1\leq k \leq n/s} q^{v_{r,n,s}(k)}
  \frac{(1-q^{-b_{r-1,n-sk}})(1-q^{-b_{r-1,k}})}{(1-q^{-1})^{2}}  \\
  & = \eta_{r,n,s} (q) \Bigl(1+ q^{-v_{r,n,s}(1)} \cdot \sum_{2\leq k
      \leq n/s} q^{v_{r,n,s}(k)}
    \frac{(1-q^{-b_{r-1,k}})(1-q^{-b_{r-1,n-sk}})}{(1-q^{-r})(1-q^{-b_{r-1,n-s}})}
  \Bigr)  \\
  & \leq \eta_{r,n,s} (q) \Bigl(1+ q^{-v_{r,n,s}(1)} \cdot \sum_{2\leq k
      \leq n/s} q^{v_{r,n,s}(k)}
    \frac{(1-q^{-b_{r-1,k}})}{(1-q^{-r})}
  \Bigr)  \label{eq:Qup} \\
& \leq \eta_{r,n,s} (q) \Bigl(1+ \frac{2q^{-v_{r,n,s}(1)+v_{r,n,s}(2)}}{(1-q^{-r})(1-q^{-1})}\Bigr) \leq \eta_{r,n,s}(q) \Bigl( 1 + \frac{16}{3}q^{-\delta}\Bigr),
\end{align}
using the bound of \autoref{lem:conv}~\ref{cor:geosum2}.

To prove \eqref{eq:32}, we observe that $\Qone{r,n,s,1}$ contains an
injective image of $(\Pone{r,n-s} \setminus \Qone{r,n-s,s})\times
\Ione{r,1}$ by $(g,h) \mapsto g\cdot h^{s}$.  For $n \geq 3s$, we get from $\Ione{r,1} = \Pone{r,1}$
\begin{align}
  \# \Qone{r,n,s}  & \geq \# \Qone{r,n,s,1}  \\
  & \geq \# \Ione{r,1} \cdot \# (\Pone{r,n-s} \setminus
  \Qone{r,n-s,s} ) \label{eq:Q1decomp}\\
  & \geq \# \Pone{r,1} \cdot (\# \Pone{r,n-s} - \# \Qone{r,n-s,s}) \\
  & \geq \eta_{r,n,s}(q) \cdot \Bigl(1 - \frac{\eta_{r,n-s,s}(q) (1+\frac{16}{3}q^{-r})}{\# \Pone{r,n-s}}\Bigr) \\
& \geq \eta_{r,n,s}(q) \cdot \Bigl(1 -
q^{b_{r,n-2s}-b_{r,n-s}+r}\frac{(1-q^{-r})(1-q^{-b_{r-1,n-2s}})(1+\frac{16}{3}q^{-r})}{(1-q^{-1})(1-q^{-b_{r-1,n-s}})}\Bigr)
\quad \label{eq:RHS} \\
& \geq \eta_{r,n,s}(q) \Bigl( 1 - \frac{7}{2} q^{-\delta-r(r-1)/2}  \Bigr), \label{eq:Qlo}
\end{align}
if $(n,s)\neq (8,2)$ using \eqref{eq:31} for $\Qone{r,n-s,s}$ with exponent $\delta \geq r$ by \autoref{thm:Q_by_gen}~\ref{item:4}.

If $(n,s) = (8,2)$, we modify \eqref{eq:RHS} according to \eqref{eq:9} and get
  \begin{align*}
    \# \Qone{r,8,2} & \geq \eta_{r,8,2}(q) \Bigl( 1 - \frac{3}{2} \bigl(1+ 2q^{-\binom{r+3}{4}-r+1}\bigr) q^{-\delta-r(r-1)/2}  \Bigr)\\
  & \geq \eta_{r,8,2}(q) ( 1 - 2q^{-\delta-r(r-1)/2} ).
  \end{align*}

Combining \eqref{eq:31} and \eqref{eq:32} proves \ref{item:9}.
\end{proof}

We note that for $(n,s)=(6,2)$, inequality \eqref{eq:45} follows from \eqref{eq:53} if $r=2$ and is false for sufficiently large $q$ if $r \geq 3$.

\begin{figure}
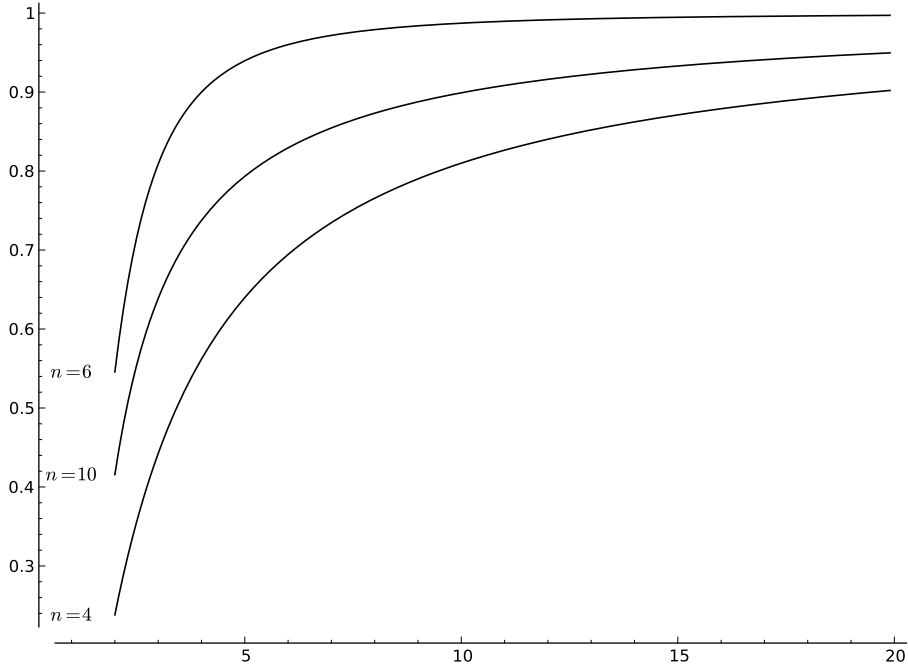

\centering
\sageplot[width=\textwidth]{QQ_vs_eta_plot(2,[4,6,10],2)}
\caption{The normalized relative error in
  \autoref{thm:Q_by_gen}~\ref{item:18}--\ref{item:13} for $(r,s)=(2,2)$.}
\label{fig:QQ_vs_eta}
\end{figure}

\autoref{fig:QQ_vs_eta} shows plots of $(\scoeff{Q}_{r,n,s}(\qvar)-\eta_{r,n,s}(\qvar))/(\eta_{r,n,s}(\qvar)
  \qvar^{-\delta})$ for $r=2$, $s=2$ and $n=4,6,10$,  as we substitute for
  $\qvar$ real numbers from $2$ to $20$.

\begin{remark}
  As noted in \autoref{rem:exp_decay} for reducible polynomials, the
  relative error term is (essentially) exponentially decreasing in the input size, and exponentially decaying in any of the parameters $r$, $n$, $s$, and $\log_2 q$, when the other three are fixed.
\end{remark}

  In the bivariate case, \citet[Theorem~3.1]{gat08-incl-gat07} approximates the quotient $\# \Qone{2,n,s}(\FF_q) / \# \Pone{2,n}(\FF_q)$ (using our notation) by
\[
q^{-(2ns-s^2+3s-4)/2}\frac{(1+q^{-1})(1-q^{-n+s-1})}{1-q^{-n-1}},
\]
which equals the term $\eta_{2,n,s}(q)/\# \Pone{2,n}(\FF_q)$ derived from our analysis above.

We append handy bounds using \autoref{cor:handy_R}.

\begin{corollary}
For $r,s,q \geq 2$, and $n \geq s$, we have
  \begin{gather}
  \frac{1}{2} q^{\binom{r+n-s}{r}+r-1} \leq \# \Qone{r,n,s} (\FF_{q}) \leq 10 
q^{\binom{r+n-s}{r}+r-1}, \\
  \frac{1}{2} q^{-\binom{r+n}{r}+\binom{r+n-s}{r}+r} \leq \frac{\# \Qone{r,n,s} (\FF_{q})}{\# \Pone{r,n} (\FF_{q})} \leq 5 
q^{-\binom{r+n}{r}+\binom{r+n-s}{r}+r}, \\
\frac{1}{6} 
q^{-\binom{r+n-1}{r}+\binom{r+n-s}{r}} \leq \frac{\# \Qone{r,n,s} (\FF_{q})}{\# \Rone{r,n} (\FF_{q})} \leq 19 
q^{-\binom{r+n-1}{r}+\binom{r+n-s}{r}}.
\end{gather}
\end{corollary}


We conclude this section with bounds for the number of $s$-powerfree polynomials.

\begin{corollary}
Let $r,s,q \geq 2$, $n \geq 0$, and $\eta_{r,n,s}$ and $\delta$ as in \autoref{thm:Q_by_gen}.  We have
\begin{equation}
  \label{eq:56}
\# \Pone{r,n}(\FF_{q})  - 3 \eta_{r,n,s}(q) \leq  \# \Sone{r,n,s} (\FF_q) \leq \# \Pone{r,n}(\FF_{q}),
\end{equation}
and more precisely
\begin{equation}
  \# \Sone{r,n,s} (\FF_q) =
\begin{cases}
\# \Pone{r,n} (\FF_{q}) & \text{for $n<s$,} \\
\# \Pone{r,n} (\FF_{q}) - \eta_{r,n,s}(q) & \text{for $s \leq n < 2s$,} \\
\# \Pone{r,n} (\FF_{q}) - \eta_{r,n,s} (q) \bigg( 1 + q^{-\delta} \cdot \frac{1-q^{-\binom{n+r-2s-1}{r-1}}}{1-q^{-\binom{n+r-s-1}{r-1}}}  & \\
\quad  \cdot \Big( \frac{1-q^{-r(r+1)/2}}{1-q^{-r}} - q^{-r(r-1)/2}\frac{1-q^{-r}}{1-q^{-1}}\Big)\bigg) &  \text{for $2s \leq n < 3s$,}
\end{cases}
\end{equation}
\begin{equation}
  \left| \# \Sone{r,6,2} (\FF_q) - (\# \Pone{r,n} (\FF_{q}) - \eta_{r,6,2}(q)) \right| \leq
  \eta_{r,6,2}(q) \cdot 2 q^{-\delta+(r-2)(r-1)(r+3)/6},
\end{equation}
and for $n \geq 3s$ with $(n,s) \neq (6,2)$
\begin{equation}
  \left| \# \Sone{r,n,s} (\FF_q) - (\# \Pone{r,n} (\FF_{q})  - \eta_{r,n,s}(q)) \right| \leq
  \eta_{r,n,s}(q) \cdot 6 q^{-\delta}.
\end{equation}
\end{corollary}


\section{Relatively irreducible polynomials}
\label{sec:rel_irr}

A polynomial over $\F$ is \emph{absolutely
  irreducible} if it is irreducible over an algebraic closure of $\F$,
and \emph{relatively irreducible} if it
is irreducible over $\F$ but factors over some extension field of $\F$.  We
define
\begin{align}
  \Aone{r,n} (\F) & = \{ f \in \Pone{r,n}(\F) \colon f \text{ is absolutely irreducible}\} \subseteq \Ione{r,n}(\F),\\
\Eone{r,n} (\F) & = \Ione{r,n} (\F) \setminus \Aone{r,n}(\F) \label{eq:66}.
\end{align}
As before, we restrict ourselves to finite fields and recall that all our polynomials are monic. For a field extension $\FF_{q^{k}}$ over $\FF_{q}$ of
degree $k$, we consider the Galois group $G_{k} =
\Gal(\FF_{q^{k}}:\FF_{q}) \cong \ZZ_k$. It acts on $\FF_{q^k}[x]$
coefficientwise and we have the ``norm'' map
\begin{align*}
  \varphi_{r,n,k} \colon \Pone{r,n/k} (\FF_{q^{k}}) & \to \Pone{r,n}
  (\FF_{q}),\\
  g & \mapsto \prod_{\sigma \in G_{k}} g^{\sigma},
\end{align*}
for each $k$ dividing $n$.  Since $(\varphi_{r,n,k} (g))^\tau = \varphi_{r,n,k}
(g)$ for any $\tau \in G_{k}$ and therefore $\varphi_{r,n,k}(g) \in \Pone{r,n}
(\FF_{q})$, this map is well-defined.

Relatively irreducible polynomials in $\Pone{r,n} (\FF_{q})$ are the product of all conjugates of an irreducible polynomial $g$ defined over some extension field $\FF_{q^{k}}$.  If $g$ itself is relatively irreducible over $\FF_{q^{k}}$, then there
 exists an appropriate multiple $j$ of $k$ and  $h \in
 \Pone{r,n/j}(\FF_{q^{j}})$ with the same image
 $\varphi_{r,n,k}(g)=\varphi_{r,n,j}(h)$ in $\Pone{r,n} (\FF_{q})$ and
 the property that $h$ is absolutely irreducible.  So, every relatively irreducible polynomial is contained in $\varphi_{r,n,k}(\Aone{r,n/k}(\FF_{q^{k}}))$ for a unique $k>1$ dividing $n$.  Furthermore, the absolutely irreducible polynomials in $\Pone{r,n}(\FF_{q})$ are exactly those in $\varphi_{r,n,1}(\Aone{r,n}(\FF_{q}))$, and we summarize
\begin{align}
\Aone{r,n}(\FF_{q}) & = \varphi_{r,n,1}(\Aone{r,n}(\FF_{q})), \label{eq:15} \\
    \Eone{r,n}(\FF_{q}) & \subseteq \bigcup_{1 < k \,\mid\, n} \varphi_{r,n,k} (\Aone{r,n/k}(\FF_{q^{k}}))   \label{eq:6}.
\end{align}
In order to replace the latter by an equality, we let
\begin{equation}
  \label{eq:3}
  \Aonep{r, n/k} (\FF_{q^{k}}) = \Aone{r,n/k} (\FF_{q^{k}}) \setminus \bigcup_{s \,\mid\, k,\,  s \neq k} \Aone{r, n/k} (\FF_{q^{s}})
\end{equation}
be the set of absolutely irreducible polynomials over $\FF_{q^{k}}$ that are not defined over a proper subfield containing $\FF_{q}$, and
\begin{equation}
  \label{eq:58}
   \Ione{r,n,k} (\FF_{q}) = \varphi_{r,n,k} (\Aonep{r,n/k}(\FF_{q^{k}})).
\end{equation}

\begin{lemma}
  \label{lem:irred}
\begin{ronumerate}
\item \label{item:6}
We have the disjoint union
\begin{equation} \label{eq:11}
  \Ione{r,n}(\FF_{q}) = \dot{\bigcup_{k \,\mid\, n}} \Ione{r,n,k}  (\FF_{q})
\end{equation}
and more precisely
\begin{align}
  \Aone{r,n} (\FF_{q}) & = \Ione{r,n,1}  (\FF_{q}), \label{eq:72} \\
\Eone{r,n} (\FF_{q}) & = \dot{\bigcup_{1 < k \,\mid\, n}} \Ione{r,n,k}  (\FF_{q}). \label{eq:union_E}
\end{align}
\item\label{item:5} \hspace*{\fill} $ \# \Ione{r,n,k}  (\FF_{q}) = \frac{1}{k} \#
  \Aonep{r,n/k}(\FF_{q^{k}})$. \hspace*{\fill}
\end{ronumerate}
\end{lemma}

\begin{proof}
  \begin{ronumerate}
\item Let $g \in \Aone{r,n/k} (\FF_{q^k})$.  By definition, $g$ is monic.  The $k$ conjugates $g^{\sigma}$, for $\sigma \in G_k$, are pairwise non-associate if and only if the coefficients are not contained
  in some proper subfield of $\FF_{q^k}$.  This shows
  \begin{equation} \label{eq:63}
  \Ione{r,n,k}  (\FF_{q})  \subseteq   \Ione{r,n}(\FF_{q}).
  \end{equation}
Let $f \in \Ione{r,n} (\FF_{q})$.  Then $f = \varphi_{r,n,k} (g)$ for some $ g\in \Aone{r,n/k} (\FF_{q^{k}})$, with $k$ dividing $n$ as observed in \eqref{eq:6}.  If $g$ has coefficients from a subfield of
   $\FF_{q^k}$, say $g \in \Aone{r,n/k} (\FF_{q^s})$ for some $s<k$
   dividing $k$, then $g^\sigma$  equals $g$ for some $\sigma \in G_k \setminus \{ \id \}$.  Taking the smallest such $s$ and
   \[
   h = \prod_{\tau \in G_s} g^{\tau} \in \Ione{r,n,k/s} (\FF_q),
   \]
   we have $h^{k/s} = \varphi_{r,n,k} (g)$. Hence
   $\varphi_{r,n,k} (g)$ is a $(k/s)$-th power and therefore reducible, in contradiction to the choice of $f$. This shows that $ g\in \Aonep{r,n/k} (\FF_{q^{k}})$ and a fortiori
\begin{equation}
\label{eq:57}
    \Ione{r,n}(\FF_{q}) \subseteq \bigcup_{k \,\mid\, n} \Ione{r,n,k}  (\FF_{q}).
  \end{equation}
The disjointness follows from the fact that the factorization of $\varphi_{r,n,k}(g)$ for any $g \in \Aonep{r,n/k}(\FF_{q^{k}})$ has exactly $k$ irreducible factors over $\FF_{q^{n}}$, and \eqref{eq:11} follows with \eqref{eq:63}.

Finally, \eqref{eq:72} and \eqref{eq:union_E} follow from \eqref{eq:15} and \eqref{eq:66}, respectively.
  \item  Let $g,h \in \Ione{r,n/k} (\FF_{q^{k}})$.  Then $\varphi_{r,n,k} (g) = \varphi_{r,n,k} (h)$ if and only if $h=g^{\sigma}$ for some automorphism $\sigma \in G_{k}$.  Sufficiency is a direct computation and necessity follows from the unique factorization of $\varphi_{r,n,k}(g)$ and $\varphi_{r,n,k}(h)$ over $\FF_{q^{k}}$.  Therefore, the
size of each fibre of $\varphi_{r,n,k}$ on $\Aonep{r,n/k}
(\FF_{q^{k}})$ is $\# G_{k} = k$. \qedhere
\end{ronumerate}
\end{proof}

We omit the parameter $r$ from the notation of the generating functions and their coefficients.  The generating function $\nAp (\FF_{q^{k}})$ of $\#
\Aonep{r,n}(\FF_{q^{k}})$ is related to the generating function
$\nA(\FF_{q})$ of $\# \Aone{r,n} (\FF_{q})$ by definition \eqref{eq:3} and we
find by inclusion-exclusion
\begin{equation}
  \label{eq:13}
  \nAp (\FF_{q^{k}})  = \sum_{s \,\mid\, k} \mu(k/s) \nA(\FF_{q^{s}}).
\end{equation}
With \eqref{eq:11} and \autoref{lem:irred}~\ref{item:5}, we relate this to the generating function $\nI (\FF_{q})$ of
irreducible polynomials as introduced in \autoref{sec:gen} and obtain
\begin{align}
  [z^{n}] \nI (\FF_{q}) & = \sum_{k \,\mid\, n } \frac{1}{k} \sum_{s \,\mid\, k} \mu (k/s) \cdot [z^{n/k}] \nA (\FF_{q^{s}}), \\
  [z^{n}] \nA (\FF_{q}) & = \sum_{k \,\mid\, n} \frac{1}{k} \sum_{s \,\mid\, k} \mu (s)
\cdot  [z^{n/k}] \nI (\FF_{q^{s}}) \label{eq:43}
\end{align}
with  M\"{o}bius inversion. 

A Maple program to compute the latter is shown in \autoref{fig:algoabs}.
\begin{figure}
\centering
\begin{Verbatim}[frame=single]
absirreds:=proc(n,r) local k,s: option remember:
    add(1/k*add(mobius(s)*subs(q=q^s,coeff(irreduciblesGF(
        z,n/k,r),z^(n/k))),s=divisors(k)),k=divisors(n))
end:

absirredsGF:=proc(z,N,r) local k,s: option remember:
    sum('absirreds(k,r)*z^k',k=1..N)
end:

relirredsGF:=proc(z,N,r) option remember:
    irreduciblesGF(z,N,r)-absirredsGF(z,N,r);
end:

relirreds:=proc(n,r)
    coeff(sort(expand(relirredsGF(z,n,r))),z^n):
end:
\end{Verbatim}
\caption{Maple program to compute the number of relatively irreducible polynomials in $r$ variables of degree $n$.}
\label{fig:algoabs}
\end{figure}
Exact values for $\# \Eone{2,n} (\FF_{q})$ with $n\leq 6$ are given in \citet[Table~4.1]{gat08-incl-gat07}.  We extend this in \autoref{tablaE}.

\begin{table}\footnotesize
\centering
\sagestr{E_table(2,8)}

\sagestr{E_table(3,7)}

\sagestr{E_table(4,6)}
  \caption{Exact values of $\# \Eone{r,n} (\FF_{q})$ for small values
    of $r$ and $n$.}
  \label{tablaE}
\end{table}

For an explicit formula, we combine the expression for $\nI_{n} (\FF_{q}) = \nI_{n}$ from \autoref{pro:R_exact_by_recursion} with \eqref{eq:43}.
\begin{theorem}
\label{pro:E_exact_by_recursion}
For $r,n \geq 1$, $q \geq 2$, $M_{n}$ as in \eqref{eq:10}, and $\nP_{n}(\FF_{q}) = \nP_{n}$ as in \eqref{eq:44}, we have
  \begin{align}
  \nA_{0} (\FF_{q})&   = 0, \\
\nA_{n} (\FF_{q})
& = - \sum_{s \, \mid \, k \,\mid\, n} \frac{\mu (s)}{k} \sum_{m \,\mid\, n/k } \frac{\mu (m)}{m} \sum_{j \in M_{n/(km)}} \frac{(-1)^{\abs{j}}}{\abs{j}} \nP_{j_{1}} (\FF_{q^{s}}) \nP_{j_{2}} (\FF_{q^{s}}) \cdots \nP_{j_{\abs{j}}} (\FF_{q^{s}}), \\
\nE_{0} (\FF_{q})& = 0, \\
\nE_{n} (\FF_{q}) & =  -\sum_{1 < k \,\mid\, n} \frac{1}{k} \sum_{s \,\mid\, k} \mu (s) \nI_{n/k} (\FF_{q^{s}}) \label{eq:60} \\
& = \sum_{1< k \,\mid\, n} \frac{1}{k} \sum_{\substack{s \,\mid\, k \\ m \,\mid\, n/k }} \frac{\mu (s) \mu (m)}{m} \sum_{j \in M_{n/(km)}} \frac{(-1)^{\abs{j}}}{\abs{j}} \nP_{j_{1}} (\FF_{q^{s}}) \nP_{j_{2}} (\FF_{q^{s}}) \cdots \nP_{j_{\abs{j}}} (\FF_{q^{s}}).
  \end{align}
\end{theorem}

We check that for $r=1$ we obtain the expected result
\begin{equation}
\nA_{n} (\FF_{q}) =
\begin{cases}
q & \text{if $n =1$}, \\
0 & \text{if $n > 1$}.
\end{cases}
\end{equation}
To this end, we use the well-known fact that
\begin{equation}
\sum_{s \,\mid\, n} \mu (s) =
\begin{cases}
1 & \text{if $n =1$}, \\
0 & \text{if $n > 1$}.
\end{cases}
\end{equation}
From \eqref{eq:14} and \eqref{eq:80} we have
\begin{align}
n \nA_{n} (\FF_{q}) & = \sum_{\substack{s \,\mid\, k \,\mid\, n \\ t \,\mid\, n/k}} \mu(s) \mu (t) q^{\frac{ns}{kt}}  = \sum_{\substack{s \,\mid\, k \,\mid\, n \\ a \,\mid\, n/k}} \mu(s) \mu (n/(ka)) q^{sa} \\
& = \smashoperator{\sum_{m \,\mid\, n}} q^{m} \smashoperator{\sum_{\substack{s \,\mid\, k \,\mid\,  n \\ m = sa,\, a \,\mid\, n/k}}} \mu(s) \mu (n/(ka))  = \smashoperator{\sum_{m \,\mid\, n}} q^{m} \smashoperator{\sum_{s\,\mid\,m}} \mu(s) \smashoperator{\sum_{\substack{s \,\mid\, k \,\mid\,  n \\ m/s \,\mid\, n/k}}} \mu (ns/(mk)) \\
& = \smashoperator{\sum_{m \,\mid\, n}} q^{m} \smashoperator{\sum_{s\,\mid\,m}} \mu(s) \smashoperator{\sum_{j \,\mid\, n/m}} \mu (n/(mj))  \\
& = \smashoperator{\sum_{m \,\mid\, n}} q^{m} \smashoperator{\sum_{s\,\mid\,m}} \mu(s) \smashoperator{\sum_{i \,\mid\, n/m}} \mu (i) =
\begin{cases}
q & \text{if $n =1$}, \\
0 & \text{if $n > 1$},
\end{cases}
\end{align}
where $a=n/(kt)$, $m=as$, $j=k/s$, and $i=n/(mj)$.

\begin{figure}
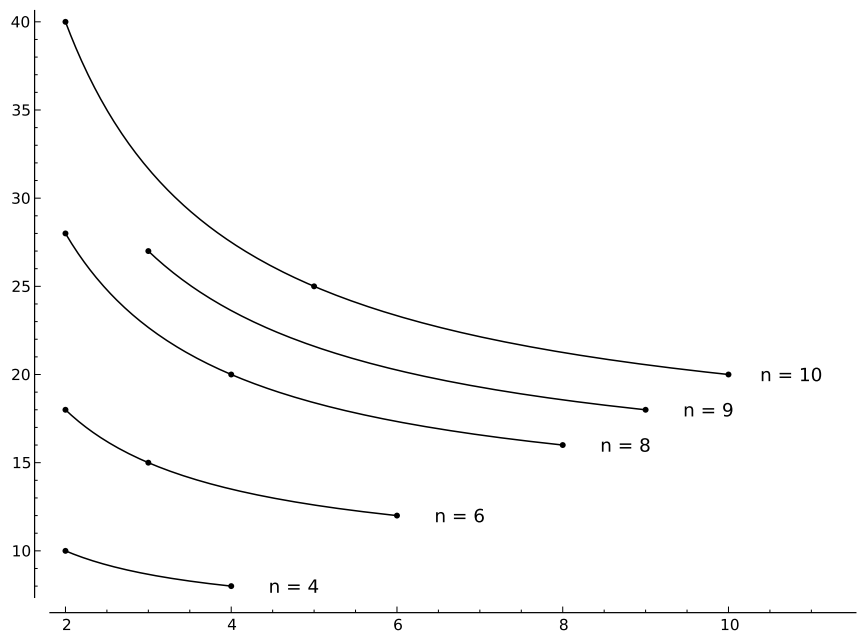

\centering
\sageplot[width=\textwidth]{w_plot(2,srange(4,11))}
\caption{Graphs for $w_{2,n}(k)$ on $[\ell,n]$ for composite $n$ in the range from $4$ to $10$, where $\ell$ denotes the smallest prime divisor of $n$.  The dots represent the values at divisors of $n$.}
\label{fig:w}
\end{figure}

The remainder of this section deals with the case $r \geq 2$.  For the approach by symbolic generating functions, we
define, with $\sI (\qvar, z)$ as in \eqref{eq:relationIP}, the two power series $\sA, \sE \in \QQ (\qvar) \bbracket{z}$ by
\begin{align}
    \scA_{0} (\qvar) & = \scI_{0} (\qvar) = 0, \\
    \scA_{n} (\qvar) & = \sum_{k \,\mid\, n} \frac{1}{k} \sum_{s \,\mid\, k} \mu (s) \sI_{n/k} (\qvar^{s}) \in \ZZ [\qvar] \text{ for } n >0, \label{eq:14} \\
\sA (\qvar, z) & = \sum_{n\geq 0} \scA_{n} (\qvar) z^{n} \in \ZZ [\qvar]\bbracket{z}, \\
\sE (\qvar, z) & = \sI (\qvar, z)- \sA(\qvar, z) \\
& = - \sum_{1 < k \,\mid\, n} \frac{1}{k} \sum_{s \,\mid\, k} \mu (s) \sI_{n/k} (\qvar^{s}) \in \ZZ [\qvar]\bbracket{z}. \label{eq:90}
\end{align}
Then
\begin{equation} \label{eq:99}
\begin{split}
  \sA_{n}(q) & = \# \Aone{r,n} (\FF_{q}), \\
  \sE_{n}(q) & = \# \Eone{r,n} (\FF_{q}).
\end{split}
\end{equation}
The inner sum of \eqref{eq:90} has degree $\degq \sI_{n/k} (\qvar^{k})$ in $\qvar$.  Let $n$ be composite and $\ell$ its smallest prime divisor.  For $k=\ell$, this inner sum consists of only two terms and we find
\begin{align}
\sE_{n} (\qvar) & = \frac{1}{\ell} (\sI_{n/\ell} ( \qvar^{\ell}) - \sI_{n/\ell}(\qvar)) - \sum_{\ell < k \,\mid\, n} \frac{1}{k} \sum_{s | k} \mu(s) \sI_{n/k} (\qvar^{s}) \\
& = \frac{1}{\ell} (\sP_{n/\ell}
( \qvar^{\ell}) - \sR_{n/\ell} ( \qvar^{\ell}) - \sI_{n/\ell}(\qvar)) + O (\qvar^{\max_{\ell < k \mid n} w_{r,n}(k)}), \label{eq:77}
\end{align}
with
\begin{equation}
 \label{eq:59} w_{r,n}(k) = \degq (\scoeff{I_{n/k}(\qvar^k)}) = \degq (\scP_{n/k} (\qvar^{k})) = k ((r+n/k)\ttf{r}/r! - 1)
\end{equation}
for any divisor $k$ of $n$.   \autoref{tab:terms_of_E} lists the degree in $\qvar$ for all summands in \eqref{eq:77}. We  consider $w_{r,n}$ as a function on the real interval $[1,n]$, see \autoref{fig:w}.

\begin{table}\footnotesize
  \centering
  \begin{tabular}{|L|L|}
\hline
    \text{summand} & \degq \\
\hline
\sP_{n/\ell} (\qvar^{\ell}) & \ell(b_{r,n/\ell}-1) =  w_{r,n} (\ell) \\
\sR_{n/\ell} ( \qvar^{\ell}) & \ell ( b_{r,n/\ell-1}+r-1) =  w_{r,n} (\ell) - \ell(b_{r-1,n/\ell}-r)\\
\sI_{n/\ell}(\qvar) & b_{r,n/\ell} -1 = \frac{1}{\ell} w_{r,n} (\ell) \\
\sum_{\ell < k \,\mid\, n} \sI_{n/k} (\qvar^{k}) & \leq \max_{\ell < k \,\mid\, n} w_{r,n}(k)  \\
\hline
  \end{tabular}
\caption{Summands of $\scE$ and their degrees in $\qvar$.}
\label{tab:terms_of_E}
\end{table}

\begin{lemma}
  \label{lem:w_mon} Let $r\geq2$, $n$ be composite, $\ell$ the
smallest and $k_{2}$ the second smallest divisor of $n$
greater than $1$.

\begin{ronumerate}
\item \label{lem:omega_inc} The function $w_{r,n}(k)$ is strictly
decreasing in $k$ on $[1, n]$.
\item \label{item:15} For composite $n \neq 4, 6$, we have
  \begin{equation} \label{eq:67}
w_{r,n}(\ell) - w_{r,n}(k_{2}) - w_{r-1,n}(\ell) \geq 0.
  \end{equation}
\item \label{item:14} For composite $n > \ell k_{2}$ different from $12$, we have
  \begin{equation}
    \label{eq:68} w_{r,n}(\ell) - w_{r,n}(k_{2}) - w_{r-1,n}(\ell) \geq \log_{2} n - 2.
  \end{equation}
This also holds if $n = 12$ and $r \geq 3$.
\end{ronumerate}
\end{lemma}

The inequality \eqref{eq:67} is false when $n$ is $4$ or $6$, and \eqref{eq:68} is false for $n=12$, $r=2$.

\begin{proof}

\ref{lem:omega_inc} We compute
  \begin{align}
 w_{r,n}'(k) &=
    \frac{(r+n/k)\ttf{r}}{r!} - \frac{n}{r! k} \sum_{1\leq i\leq r} \frac{(r+n/k)\ttf{r}}{i+ n/k} - 1  \\
    & = \frac{(r+n/k)\ttf{r}}{r!} \Bigl( 1 - \sum_{1\leq i\leq r}
      \frac{1}{1+i\frac{k}{n}} \Bigr) -1. \label{eq:w_prime}
  \end{align}

  If $r \geq 3$, then
  \[
  \sum_{1 \leq i\leq r} \frac{1}{1+i\frac{k}{n}} \geq \sum_{1\leq
    i\leq 3} \frac{1}{1+i} > 1
  \]
for all $1 \leq k \leq  n$, which proves $w_{r,n}'(k) < 0$.

  If $r = 2$, we evaluate \eqref{eq:w_prime} as
  \begin{equation}
w_{2,n}'(k) = \frac{(1+n/k)(2+n/k)}{2}\Bigl(1-\frac{1}{1+k/n}-\frac{1}{1+2k/n}\Bigr)-1
 = -\frac{n^{2}}{2k^{2}}
  \end{equation}
to find $w_{2,n}'(k) < 0$ for all $k$.

For \ref{item:15} and \ref{item:14}, we first show that the sequence $a_{r,n} =  w_{r,n}(\ell) - w_{r-1,n}(\ell) - w_{r,n}(k_{2}) = \ell b_{r,n/l-1} -k_{2} (b_{r,n/k_{2}}-1)$ is monotonically increasing in $r$.  We have
  \begin{equation}
   a_{r,n}-a_{r-1,n} = \ell b_{r,n/\ell-2} - k_{2} b_{r,n/k_{2}-1} \geq 0
  \end{equation}
if and only if
  \begin{equation}
    \label{eq:quotient}
A_{r,n} = \frac{\ell(r+n/\ell - 2)\ttf{r}}{k_{2} (r+n/k_{2} - 1)\ttf{r}} \geq 1
  \end{equation}
and prove the latter by induction on $r \geq 2$.

For $r=2$, we have to prove
\begin{equation}
n(k_{2} -\ell ) \geq 2 \ell k_{2}. \label{eq:omega_increasing}
\end{equation}
If $k_{2} = \ell+1$, then $\ell = 2$, $k_{2} =3$ and since we exclude
$n=6$, we have $n \geq 12$ to show \eqref{eq:omega_increasing}.  If
$k_{2} \geq \ell+2$, we distinguish two cases.  Now, $k_{2} = n$ if
and only if $n = \ell^2$.  Since we exclude $n=4$, we then have $\ell
\geq 3$ and \eqref{eq:omega_increasing} follows.  If $k_{2} \neq n$,
then $k_{2} \leq \sqrt{n} < n$ and therefore $2 \ell k_{2} < 2
\sqrt{n} \sqrt{n} \leq (k_{2} - \ell) n $.

For the induction step, we have
  \[
  A_{r,n} = A_{r-1,n} \frac{n/\ell -2 +r}{n/k_{2} - 1+r} \geq  \frac{n/\ell -2 +r}{n/k_{2} - 1+r} \geq 1,
  \]
where the last inequality is equivalent to $n(k_{2} -\ell ) \geq \ell
k_{2}$, which follows from \eqref{eq:omega_increasing}.

With this monotonicity of $a_{r,n}$ in $r$, it is sufficient to check
\ref{item:15} and \ref{item:14} for the smallest admissible value of $r$.

\begin{ronumerate}
\item[\ref{item:15}] We have
  \begin{equation}
  \label{eq:78}
    a_{2,n} = \frac{n}{2} \Bigl(\frac{n}{\ell}-\frac{n}{k_{2}}-2\Bigr).
  \end{equation}
For
\begin{itemize}
\item  $n =\ell^{2}$, $\ell \neq 4$,
\item $n=\ell k_{2}$, $n \neq 6$, or
\item $n=12$,
\end{itemize}
this is non-negative by direct computation, and in the remaining case, $n >
\ell k_{2}$ different from $12$, by \ref{item:14}.
\item[\ref{item:14}] For $n > \ell k_{2}$ different from $12$, we have
  $n/\ell - n/k_{2} \geq 3$ and find with \eqref{eq:78}
  \begin{equation}
    \label{eq:110}
     a_{2,n} \geq \frac{n}{2} > \log_{2}n -2.
  \end{equation}
For $n=12$ and $r \geq 3$, we compute directly $a_{3,12} = 10 >
\log_{2} 12 - 2$. \qedhere
  \end{ronumerate}
\end{proof}

This lemma allows us to order the summands in \eqref{eq:77} by $\degq$, and the approach by generating functions gives the following result.

\begin{theorem}
\label{thm:E_by_gen}
Let $r, n \geq 2$, let $\ell$ be the smallest prime divisor of $n$, and
 \begin{align}
    \epsilon_{r,n}(\qvar) & = \frac{\qvar^{\ell ( \binom{r+n/\ell}{r} - 1)}}{\ell(1-\qvar^{-\ell})} \in \QQ(\qvar), \label{eq:epsilon} \\
\kappa & = (\ell-1) (\binom{r-1+n/\ell}{r-1}-r ) + 1.
 \end{align}
Then the following hold.
\begin{ronumerate}
  \item \label{enum3:i} $\scoeff{E_{1}} (\qvar)= 0$.
  \item \label{enum3:ii} If $n$ is prime, then
    \begin{align}
    \scoeff{E_n} (\qvar) & = \epsilon_{r,n}(\qvar) (1-\qvar^{-nr}) \Bigl( 1 -
      \qvar^{-r(n-1)}\frac{(1-\qvar^{-r})(1-\qvar^{-n})}{(1-\qvar^{-1})(1-\qvar^{-nr})}\Bigr).
    \end{align}
    \item \label{enum3:iii} If $n$ is composite, then $\kappa \geq 2$ and
\[
    \scoeff{E_n} (\qvar) = \epsilon_{r,n}(\qvar) (
 1 + O (\qvar^{- \kappa}) ).
\]
\end{ronumerate}
\end{theorem}

\begin{proof}
For $n=1$, the sum \eqref{eq:60} is empty and this shows \ref{enum3:i}. For $n=\ell$ prime, \eqref{eq:60} simplifies to
$\sE_{n} (\qvar) = (\sI_{1} ( \qvar^{\ell}) - \sI_{1}(\qvar))/\ell  = (\sP_{1}
( \qvar^{\ell}) - \sP_{1}(\qvar))/\ell $, since $\sI_{1}=\sP_{1}$ by \autoref{thm:R-from-gen} and \ref{enum3:ii} follows.

For composite $n$, the product $(\ell-1)(b_{r-1,n/\ell}-r)$ is positive and therefore $\kappa \geq 2$.  We recall the summands of \eqref{eq:77} in \autoref{tab:terms_of_E}.  \autoref{lem:w_mon} \ref{lem:omega_inc} shows that $\max_{\ell < k \,\mid\, n} w_{r,n} (k) = w_{r,n} (k_{2})$ and we find
\begin{align}
\sE_{n} (\qvar) & = \frac{1}{\ell} (\sP_{n/\ell}
( \qvar^{\ell}) - \sR_{n/\ell} ( \qvar^{\ell}) - \sI_{n/\ell}(\qvar)) + O (\qvar^{w_{r,n}(k_{2})}).
\end{align}
Since $b_{r-1,n/\ell}-r > 0$ for composite $n$, we identify with \autoref{lem:w_mon} \ref{lem:omega_inc} as main term $\sP_{n/\ell} (\qvar^{\ell})/\ell = \epsilon_{r,n}(\qvar) (1-\qvar^{-\ell b_{r-1,n/\ell}})$.  For the summands of
\begin{equation}
  \label{eq:4}
  \sE_{n} (\qvar) / \epsilon_{r,n}(\qvar) = (1-\qvar^{-\ell b_{r-1,n/\ell}}) \Bigl(1  - \frac{\sR_{n/\ell} ( \qvar^{\ell})}{\sP_{n/\ell}
( \qvar^{\ell})} - \frac{\sI_{n/\ell}(\qvar)}{\sP_{n/\ell}
( \qvar^{\ell})}\Bigr) + O (\qvar^{w_{r,n}(k_{2})-\degq \sP_{n/\ell} (\qvar^{\ell})})
\end{equation}
we find as degrees in $\qvar$
\begin{align}
-\ell b_{r-1,n/\ell} & \leq -\kappa, \\
  \degq \sR_{n/\ell} ( \qvar^{\ell}) - \degq\sP_{n/\ell} ( \qvar^{\ell})  = -\ell (b_{r-1,n/\ell}-r)  & \leq -\kappa, \label{eq:74} \\
\degq  \sI_{n/\ell}(\qvar)   - \degq\sP_{n/\ell} ( \qvar^{\ell})  = -(\ell-1) (b_{r,n/\ell}-1) & \leq -\kappa, \label{eq:75} \\
w_{r,n}(k_{2}) - \degq\sP_{n/\ell} ( \qvar^{\ell})  \leq -\ell (b_{r-1,n/\ell}-1) & \leq -\kappa  \label{eq:76}
\end{align}
for $n \neq 4,6$ by \autoref{lem:w_mon} \ref{item:15}.  When $n$ is $4$ or $6$, the last inequality in \eqref{eq:76} is false, but still
\begin{equation}
  \label{eq:73}
  w_{r,n}(k_{2}) - \degq\sP_{n/\ell} ( \qvar^{\ell}) \leq
  -\kappa. \mbox{\qedhere}
\end{equation}
\end{proof}

On closer inspection, it is possible to partition for each
composite $n$ the range for $r$ into two non-empty intervals, where either the difference in \eqref{eq:74} or the difference in \eqref{eq:75} dominates all others.  This provides tighter bounds at the cost of further case distinctions.

The combinatorial approach yields the following result.

\begin{theorem}
  \label{thm:E_complete}
  Let $r,q \geq 2$, and $\epsilon_{r,n}$ and $\kappa$ as in \autoref{thm:E_by_gen}.
  \begin{ronumerate}
  \item \label{item:11}
    $\# \Eone{r,1} (\FF_q) = 0$.
  \item \label{item:12} If $n$ is prime, then
    \begin{gather}
    \# \Eone{r,n} (\FF_q) = \epsilon_{r,n}(q) (1-q^{-nr}) \Bigl( 1 -
      q^{-r(n-1)}\frac{(1-q^{-r})(1-q^{-n})}{(1-q^{-1})(1-q^{-nr})}\Bigr),
    \qquad \label{eq:69}
    \\
    0 \leq \epsilon_{r,n}(q) - \# \Eone{r,n} (\FF_q) \leq 3
    q^{-r(n-1)}.
  \end{gather}
  \item \label{item:16} If $n$ is composite, then
    \[
    \left| \# \Eone{r,n}(\FF_q) - \epsilon_{r,n}(q) \right| \leq
    \epsilon_{r,n}(q) \cdot 3 q^{-\kappa}.
    \]
  \end{ronumerate}
\end{theorem}

\begin{proof}
The exact statements of \ref{item:11} and \ref{item:12} were already shown in \autoref{thm:E_by_gen} and in \eqref{eq:69} we have $q^{-r(n-1)}/16$ as upper bound for $q^{-nr}$ and $32q^{-r(n-1)}/15$ as upper bound for the last subtracted term.

For \ref{item:16}, let $\ell$ be the smallest and $k_{2}$ the second smallest divisor of $n$ greater than $1$. We prove that
\begin{align}
    \# \Eone{r,n}(\FF_q) & \geq \epsilon_{r,n}(q) ( 1 - 3
      q^{-\kappa} ),   \label{pro:E_lb} \\
    \# \Eone{r,n} (\FF_q) & \leq \epsilon_{r,n}(q) (1 + 2 q^{-\ell
      (b_{r-1,n/\ell}-1)})  \quad \text{for $n \neq 4, 6$}, \label{pro:E_ub} \\
    \# \Eone{r,n} (\FF_q) & \leq \epsilon_{r,n}(q) (1 + q^{-\kappa})  \quad \text{for $n = 4, 6$.} \label{eq:79}
  \end{align}

We begin with \eqref{pro:E_lb} and have from \autoref{lem:irred}~\ref{item:5}
  \begin{align}
    \# \Eone{r,n} (\FF_{q}) & \geq \# \Ione{r,n,\ell} (\FF_{q}) = \frac{1}{\ell} \#
    \Aonep{r,n/\ell} (\FF_{q^\ell}) \\
& = \frac{1}{\ell} ( \# \Ione{r,n/\ell} (\FF_{q^{\ell}}) - \# \Ione{r,n/\ell} (\FF_{q})),
  \end{align}
  since $\ell$ is prime and there are no proper intermediate fields
  between $\FF_{q}$ and $\FF_{q^{\ell}}$.  With the
  lower bound on the number of irreducible polynomials from
  \autoref{cor:irred} this yields
  \begin{align}
    \# \Eone{r,n} (\FF_{q}) & \geq \frac{1}{\ell} ( \# \Pone{r,n/\ell}(\FF_{q^\ell}) - 2 \rho_{r,n/\ell}(q^\ell) - \# \Pone{r,n/\ell}(\FF_q) )  \\
    & = \epsilon_{r,n}(q) \Bigl( 1 - q^{-\ell b_{r-1,n/\ell}}
    - 2 q^{-\ell(b_{r-1,n/\ell}-r)} \frac{1-q^{-\ell r}}{1-q^{-\ell}}   \\
    & \quad  - q^{-(\ell -1)(b_{r,n/\ell}-1)} \frac{(1-q^{-b_{r-1,n/\ell}})(1-q^{-\ell})}{1-q^{-1}} \Bigr)  \\
    & = \epsilon_{r,n}(q) \biggl( 1 - q^{-\kappa}
    \Big(q^{-b_{r-1,n/\ell}-\ell r + 1}
+ 2 q^{-b_{r-1,n/\ell}+r+1} \frac{1-q^{-\ell r}}{1-q^{-\ell}} \\
& \quad + q^{-(\ell -1)b_{r,n/\ell-1}-\ell r + \ell +r}      \frac{(1-q^{-b_{r-1,n/\ell}})(1-q^{-\ell})}{1-q^{-1}} \Big)\biggr) \\
& \geq \epsilon_{r,n}(q) (1- q^{-\kappa} (1/16 + 8/3 + 1/4)) \\
& \geq \epsilon_{r,n}(q) (1 - 3q^{-\kappa} ).
\end{align}

For the lower bounds \eqref{pro:E_ub} and \eqref{eq:79}, we have from \autoref{lem:irred}~\ref{item:5}
\begin{align}
  \# \Ione{r,n,k}  (\FF_{q}) & = \frac{1}{k} \#  \Aonep{r,n/k} (\FF_{q^{k}})  \\
  & \leq \frac{1}{k} \#  \Pone{r,n/k} (\FF_{q^{k}})  \\
  & = q^{w_{r,n}(k)}
  \frac{1-q^{-k\binom{n/k+r-1}{r-1}}}{k(1-q^{-k})}, \label{eq:est_E_rnk}
\end{align}
with $w_{r,n}(k)$ as defined in \eqref{eq:59}.
We obtain with \eqref{eq:union_E}
  \begin{align}
    \# \Eone{r,n} (\FF_{q}) & \leq \sum_{1< k|n} \# \Ione{r,n,k}  (\FF_{q})  \\
    & \leq \sum_{1< k \,\mid\, n} q^{w_{r,n}(k)} \cdot \frac{1-q^{-kb_{r-1,n/k}}}{k(1-q^{-k})}  \\
    & =  q^{w_{r,n}(\ell)} \frac{1-q^{-\ell b_{r-1,n/\ell}}}{\ell(1-q^{-\ell})}+ \sum_{\ell < k \,\mid\, n} q^{w_{r,n}(k)}\frac{1-q^{-kb_{r-1,n/k}}}{k(1-q^{-k})}  \\
    & = \epsilon_{r,n}(q) (1-q^{-\ell b_{r-1,n/\ell}}) \\
& \quad \cdot \Bigl(1 + q^{-w_{r,n}(\ell)}
      \sum_{\ell < k \,\mid\, n} q^{w_{r,n}(k)} \frac{\ell (1 - q^{-\ell})
        (1-q^{-kb_{r-1,n/k}})} {k(1-q^{-k})(1-q^{-\ell
          b_{r-1,n/\ell}})}
    \Bigr)  \\
    & \leq \epsilon_{r,n}(q) \Bigl(1 + q^{-w_{r,n}(\ell)} \sum_{\ell <
        k \,\mid\, n} \frac{\ell}{k} q^{w_{r,n}(k)}
    \Bigr), \label{eq:64}
  \end{align}
  since $(1-q^{-k})/(1-q^{-kb_{r-1,n/k}})$ is monotone increasing with
  $k$. 

For $n=\ell^{2}$ or $n=\ell k_{2}$, we compute directly from \eqref{eq:64}
  \begin{align}
    \# \Eone{r,\ell^{2}} (\FF_q) & \leq \epsilon_{r,n}(q) \Bigl(1 +
      \frac{1}{\ell} q^{-w_{r,n}(\ell) + w_{r,n}(n)}\Bigr), \\
    \# \Eone{r,\ell k_{2}} (\FF_q)     & \leq \epsilon_{r,n}(q) \Bigl(1 + q^{-w_{r,n}(\ell) + w_{r,n}(k_{2})} \Bigl(\frac{\ell}{k_{2}} + \frac{\ell}{n}\Bigr)\Bigr) \\
    & \leq  \epsilon_{r,n}(q) (1 + q^{-w_{r,n}(\ell) + w_{r,n}(k_{2})} ),
  \end{align}
respectively.  These prove \eqref{pro:E_ub} for $n \neq 4, 6$, since
$-w_{r,n}(\ell) + w_{r,n}(k_{2}) \leq - w_{r-1,n}(\ell) \leq -\kappa$
by \autoref{lem:w_mon} \ref{item:15}, and they also show \eqref{eq:79}
for $n=4,6$ with \eqref{eq:73}.

For $n > \ell k_{2}$, we show
  \begin{equation}
 \label{eq:100}
    q^{-w_{r,n}(\ell)} \sum_{\ell<k \,\mid\, n} \frac{\ell}{k} q^{w_{r,n}(k)} \leq 2q^{-w_{r-1,n}(\ell)}.
  \end{equation}
  We use the coarse bound $\# \{k \colon \ell < k \,\mid\, n\} \leq n/2 =
      2^{\log_{2} n - 1 } \leq 2q^{\log_{2} n - 2}$ and show the stronger
\begin{align}
    q^{-w_{r,n}(\ell)} 2q^{\log_{2} n - 2} q^{w_{r,n}(k_{2})} & \leq 2q^{-w_{r-1,n}(\ell)} \\
    \intertext{or equivalently}
    -w_{r,n}(\ell) + w_{r,n}(k_{2}) & \leq -w_{r-1,n}(\ell) - \log_{2} n + 2. \label{eq:key_estimate_exp}
\end{align}
For $n \neq 12$ or $n=12$ and $r \geq 3$, this follows from
\autoref{lem:w_mon} \ref{item:14}.  For $r=2$ and $n = 12$, it suffices to evaluate left-
and right-hand side of \eqref{eq:100} to find $5/6q^{-12}<2q^{-12}$ as claimed.

Finally, we combine the bounds \eqref{pro:E_lb}, \eqref{pro:E_ub}, and \eqref{eq:79} with $-w_{r-1,n}(\ell) \leq -\kappa$ from \eqref{eq:76}.
\end{proof}

\begin{figure}
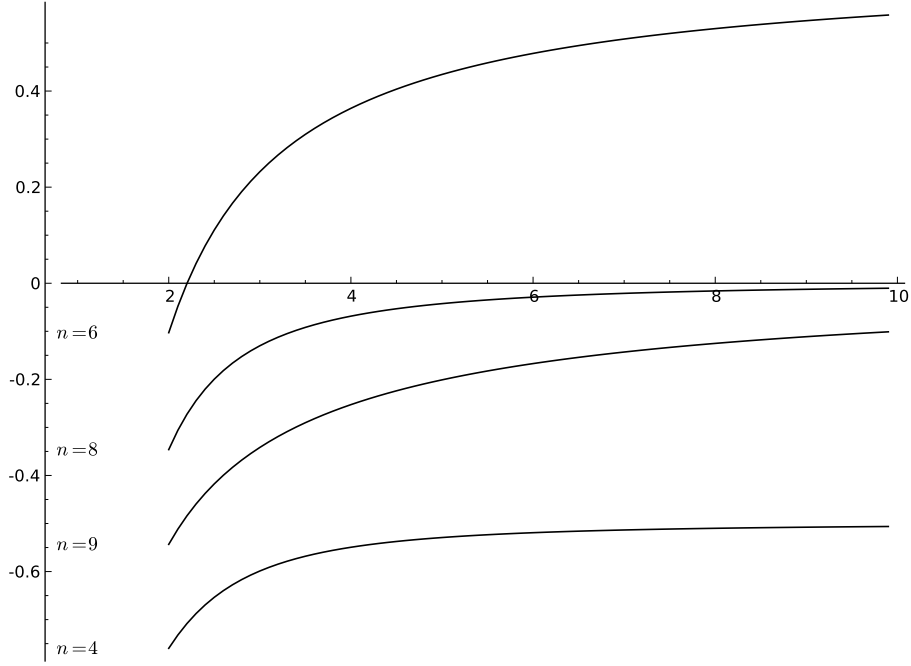

\centering
\sageplot[width=\textwidth]{EE_vs_eps_plot(2,[4,6,8,9])}
\caption{The normalized relative error in
  \autoref{thm:E_by_gen}~\ref{enum3:iii} for $r=2$.}
\label{fig:EE_vs_eps}
\end{figure}

\autoref{fig:EE_vs_eps} shows plots of $(\scoeff{E}_{r,n}(\qvar)-\epsilon_{r,n}(\qvar))/(\epsilon_{r,n}(\qvar)
  \qvar^{-\kappa})$ for $r=2$ and $n=4,6,8,9$, as we substitute for
  $\qvar$ real numbers from $2$ to $10$.

\begin{remark}
 The bivariate result of \cite{gat08-incl-gat07} approximates the ratio $\# \Eone{2,n}(\FF_q)/ \# \Pone{2,n} (\FF_q)$ by
\[
\frac{q^{-n^2(\ell - 1)/(2\ell)} (1-q^{-1})}{\ell (1-q^{-\ell})(1-q^{-n-1})}.
\]
This differs from the approximation by $\epsilon_{2,n}(q) / \# \Pone{2,n} (\FF_q)$ in \autoref{thm:E_complete} by a factor of $1-q^{-n-1}$.
\end{remark}

We append some handy bounds.
\begin{corollary}
  Let $r,n,q\geq 2$, and $\ell$ be the smallest prime divisor of $n$, then
  \begin{gather}
    \frac{1}{4 \ell} q^{\ell \binom{r+ n/\ell}{r}-\ell}\leq \# \Eone{r,n}(\FF_{q}) \leq \frac{2}{\ell} 
q^{\ell \binom{l + n/\ell}{r}-\ell}, \\
    \frac{1}{8 \ell} q^{- \binom{r+n}{r} + \ell \binom{r + n/\ell}{r} -\ell + 1}\leq \frac{\# \Eone{r,n}(\FF_{q})}{\# \Pone{r,n}(\FF_{q})} \leq  \frac{2}{\ell} 
q^{- \binom{r+n}{r} + \ell \binom{r + n/\ell}{r} - \ell + 1}, \\
    \frac{1}{8 \ell} q^{- \binom{r+n}{r} + \ell \binom{r + n/\ell}{r} - \ell + 1} \leq \frac{\# \Eone{r,n} (\FF_{q})}{\# \Ione{r,n} (\FF_{q})} \leq  \frac{2}{\ell} 
q^{- \binom{r+n}{r} + \ell \binom{r + n/\ell}{r} - \ell + 1}.
  \end{gather}
\end{corollary}


The last inequalities follow with
\autoref{cor:handy_R} for $n \geq 5$ and by computation with the exact
expressions otherwise.

We conclude with bounds for the number of absolutely irreducible polynomials by combining \autoref{cor:irred} and \autoref{thm:E_complete}.

\begin{corollary}
Let $r,n,q \geq 2$, and $\rho_{r,n}(q)$ as in \eqref{eq:rho}.  Then
\[
\# \Pone{r,n}(\FF_q) - 4 \rho_{r,n}(q) \leq \# \Aone{r,n} (\FF_{q})
\leq \# \Ione{r,n} (\FF_{q}) \leq \# \Pone{r,n} (\FF_{q}),
\]
where the 4 can be replaced by 3 for $n \geq 3$.
\end{corollary}


\section{Acknowledgments}

Joachim von zur Gathen and Alfredo Viola thank the late Philippe Flajolet for useful
discussions about Bender's method in analytic combinatorics of
divergent series in April~2008.  The work of Joachim von zur Gathen
and Konstantin Ziegler was supported by the B-IT Foundation and the
Land Nordrhein-Westfalen.  We thank the anonymous referees for their useful comments.

\bibliography{journals,refs,lncs}

\bibliographystyle{cc2}

\end{document}